\documentclass[final]{article}

\usepackage[a4paper,top=4.0cm,bottom=2.0cm,left=2.0cm,right=2.0cm]{geometry}
\usepackage{appendix}
\usepackage{pgfplots}
\usepackage{amsmath,amsfonts,mathrsfs,amsfonts,ulem}
\usepgfplotslibrary{patchplots}
\pgfplotsset{compat=1.15}
\usepackage{lineno}
\usepackage{array}
\usepackage{longtable}
\usepackage{units}
\usepackage{hyperref}
\usepackage{cleveref}
\usepackage{amsfonts}
\usepackage{listings}
\usepackage{epsfig,graphicx}
\usepackage{color}

\newcommand{\R}{{\mathbb R}}
\usepackage{tikz,tikz-3dplot}
\usetikzlibrary{shadows}
\usetikzlibrary{shapes.geometric, arrows.meta}
\usetikzlibrary{positioning}
\usetikzlibrary{calc,shapes.multipart,chains}
\usetikzlibrary{patterns.meta}
\usepackage{svg}
\usepackage{datetime}

\usepackage{algorithm}
\usepackage{algpseudocode}

\usepackage{color, bbold,bm,blindtext}
\definecolor{darkblue}{RGB}{0,125,175}
\definecolor{cherryred}{RGB}{204, 0, 51}
\definecolor{darkgreen}{RGB}{0, 204, 102}
\definecolor{darkorange}{RGB}{255, 140, 0}

\usepackage{authblk}

\allowdisplaybreaks

\begin{document}
\title{Adaptive tempering schedules with approximative intermediate measures for filtering problems}

\author[1]{Iris Rammelmüller}
\author[2]{Gottfried Hastermann}
\author[3]{Jana de Wiljes}
\affil[1]{\footnotesize Department of Statistics, University of Klagenfurt, Universitätsstraße 65–67, 9020 Klagenfurt, Austria,  \url{iris.rammelmueller@aau.at}}
\affil[2]{Institute for Mathematics, Ilmenau University of Technology, Weimarer Str. 25, 98693 Ilmenau, \url{gottfried.hastermann@tu-ilmenau.de}}
\affil[3]{\footnotesize Institute for Mathematics, Ilmenau University of Technology, Weimarer Str. 25, 98693 Ilmenau, \url{jana.de-wiljes@tu-ilmenau.de}, School of Engineering
Sciences, Department of Computational Engineering, Yliopistonkatu 34
53850 Lappeenranta, Finland LUT University, Finland}

\date{Preprint, May 2024}

\maketitle

\begin{abstract} 
Data assimilation algorithms integrate prior information from numerical model simulations with observed data. Ensemble-based filters, regarded as state-of-the-art, are widely employed for large-scale estimation tasks in disciplines such as geoscience and meteorology. Despite their inability to produce the true posterior distribution for nonlinear systems, their robustness and capacity for state tracking are noteworthy. In contrast, Particle filters yield the correct distribution in the ensemble limit but require substantially larger ensemble sizes than ensemble-based filters to maintain stability in higher-dimensional spaces. It is essential to transcend traditional Gaussian assumptions to achieve realistic quantification of uncertainties. One approach involves the hybridisation of filters, facilitated by tempering, to harness the complementary strengths of different filters. A new adaptive tempering method is proposed to tune the underlying schedule, aiming to systematically surpass the performance previously achieved. Although promising numerical results for certain filter combinations in toy examples exist in the literature, the tuning of hyperparameters presents a considerable challenge. A deeper understanding of these interactions is crucial for practical applications.

\end{abstract}
\noindent {\bf AMS(MOS) subject classifications.} 65C05, 62M20, 63G35, 
\section{Introduction}
Estimation of large-scale numerical weather and climate prediction signals on the basis of available observation and evolution models is approached with data assimilation methods, c.f. \cite{Reich2015,Law2015, Evensen2022}. These methods are typically based on the assumption that the associated posterior distribution is a Gaussian. This specification leads to highly robust filters, which possess strong tracking capabilities, c.f. \cite{deWiljes2018,deWiljes2020}. Many of the systems analysed exhibit characteristics closely aligned with near Gaussian distributions; however, they also display pronounced nonlinear behaviour, especially in case of (very) sparse and uncertain observation. Forecasting of precipitation serves as one of the most prominent and practical examples. More generally, nonlinear transport with little dissipative effects is the foundation of multiple examples with strong non-Gaussian prior distributions. Depending on the observational model the posterior might then also be substantially non-Gaussian. Considering that numerous decisions impacting public health rely on these predictions, it becomes imperative to refine approximations of the underlying uncertainties in both prior and posterior distributions of each filtering step. This refinement can be realised through non-parametric methods, which provide empirical estimates of the distributions in question. However, vastness of the state space often complicates the direct use of standard techniques such as sequential Monte Carlo methods. Consequently, substantial research has been devoted to developing methodologies that not only yield more accurate uncertainty estimates but also remain computationally viable, c.f. \cite{Beskos_Crisan_Jasra_Kamatani_Zhou_2017}. A particular line of research has focused on developing and training cost-effective surrogate models. Various paths are taken in the derivation of these models, including filtering on multiple levels of coarse-grained grids \cite{Hako2020}, estimating parameters for alternative families of functions like neural networks or systems of linear equations, c.f.\cite{KarenVeroy_etal2023, TruongVinhe2023}. These approaches are often aimed at facilitating the generation of larger ensembles of samples. To ensure that coarse-grained simulations accurately reflect the true system, noise carrying the system's physical information is incorporated into the evolution equation thereby improving the efficiency and scalability of simulations while maintaining the same accuracy in generating samples of the posterior. While this ansatz is particularly promising for specific large-scale systems, e.g. \cite{Cotter2020,Cotter2020a}, it is not trivial to adapt to other model equations and requires a series of other tools namely jittering (or rejuvenation/ inflation techniques) and tempering, which in general improve the feasibility of non-parametric empirical inference approaches. One of the key advantages of most Gaussian approximate filters opposed to standard particle filters is that individual ensemble methods are moved towards the observation in the equations, which can be interpreted as evolution equations of the particles in time in a continuous setting.  It is feasible to develop analogous evolution equations for the non-parametric scenario, leading to the introduction of filters such as the Feedback Particle Filter \cite{Taghvaei2016}, see \cite{Pathiraja2021} for an extensive and thorough generalisation. While these approaches have a similar update structure as their Gaussian counterparts they are typically limited by the existing computational boundaries. As the aim is to map samples from the prior, associated with the underlying model, to a posterior distribution conditioned on the important additional observation information source a family of methods is focused on estimating the underlying transport map \cite{spantini2022, Cheng2013,Papamakarios2021}. These maps are calculated based on the available particles and associated importance weights. However, in high dimensions, the samples may not sufficiently cover the underlying state space. Therefore, it is often necessary to either use one of the aforementioned techniques or to combine complementing filters. The latter integration is facilitated through a procedure known as tempering, which involves dividing the likelihood into parts and then merging non-parametric with parametric filters, c.f. \cite{Chustagulprom2016,Frei2013}. This blend of filters is often termed hybrid, which should not be confused with the usage of 'hybrid filters' in the context of combining ensemble and variational methods such as 3DVar or 4DVar \cite{Gharamti2021}. 
One can view these approaches in the context of preconditioners or, improved proposal densities. They have been successfully numerically implemented for specific combinations of filters for a few toy models \cite{Chustagulprom2016,Acevedo2017,Nerger2021}. 
Although these numerical experiments highlight the approach's potential, it is important to extensively examine the impacts of the ordering of the filters, fine-tuning hyperparameters such as the tempering schedule, and testing the underlying hypotheses across a broader spectrum of filters to enable application in real-life data assimilation set-ups. To explore the effects of various tempered combinations more effectively, we designed a series of toy setups to numerically investigate underlying systems. Specifically, we employ a Langevin model to simulate problems characteristic of high-dimensional spaces. We also propose a novel method for adaptively tuning the tempering schedule, which shows promising results across all toy models and filter combinations compared to previous approaches using the effective sample size.

\section{Mathematical framework} \label{sec:math}
 Let $z_k \in \R^{N_z}$ be the state of a dynamical system and $\Psi: \R^{N_z} \to \R^{N_z}$ be the solution operator of a dynamical system. Given an initial condition $z_0 \in \R^{N_z}$, a dynamical system defines a solution sequence $\lbrace z_k \rbrace_{k \geq 0}$ via the forward iteration
\begin{equation}
    z_{k+1} = \Psi(z_k).   \label{eq:evolution}
\end{equation}
\noindent
Let $\lbrace y_k \rbrace_{k \geq 0} =y \in \R^{N_y}$ be the noisy observations then at each discrete time instance we observe a possibly nonlinear function of the signal
\begin{equation}
     y_{k+1} = H(z_{k+1}) + \eta_{k+1},    \label{eq:datamodel}
\end{equation}
\noindent 
where $H:\R^{N_z} \to \R^{N_y} $ is the observation operator and $\eta_{k+1}$ is Gaussian distributed with mean zero and error covariance matrix $R\in \mathbb{R}^{N_{y} \times N_{y}}$. Given this discrete scenario, we would like to solve the time-discrete filtering problem, i.e. to estimate the density $p(z_{k} \vert y_{1:k})$ associated with the signal $z_k$ conditioned on the observation up to time $k$. Employing Bayes' theorem allows us to leverage the following relationship:
\begin{equation}\label{eq:filtering_problem}
    p(z_{k} \vert y_{1:k}) \propto p(y_{k} \vert z_{k}) p(z_{k} \vert y_{1:k-1}),
\end{equation}
where \( p(y_{k} \vert z_{k}) \) is the likelihood and \( p(z_{k} \vert y_{1:k-1}) \) is the prior distribution c.f. \cite{Jazwinski2007}. 
This representation of the posterior allows one to modify the likelihood. For instance one can introduce artificial measures 
\begin{equation}\label{eq:artificial_measures}
\frac{d\mu_r}{d\mu_0}\propto p(y_{k} \vert z_{k})^{\tau_r}
    \end{equation}
where $0=\tau_0< \ldots < \tau_r < \ldots < \tau_T=1$ is the schedule, which represents the artificial time steps in between the two actual time $k-1$ and $k$, $\mu_T=p(z_{k} \vert y_{1:k})$ and $\mu_0=p(z_{k} \vert y_{1:k-1})$. Then one can use this split to sequentially update with the intermediate measures  
\begin{equation}\label{eq:intermediate_measures}
\frac{d\mu_r}{d\mu_{r-1}}\propto p(y_{k} \vert z_{k})^{\tau_r-\tau_{r-1}}.
    \end{equation}
This procedure, known as tempering, is particularly beneficial when the supports of the prior and the posterior distributions do not overlap. Although this technique is employed broadly within the context of sequential Monte Carlo methods, c.f. \cite{Chopin2020}, it is especially crucial for filtering problems \cite{Cotter2020a,Beskos_Crisan_Jasra_Kamatani_Zhou_2017}. Here, we consider a specialised version where different approximation schemes are used for the intermediate updates. This approach shows very promising results, though the benefits vary depending on the quality of the filters used, the underlying system dynamics, and the available data \cite{Chustagulprom2016,Frei2013}. Therefore, for these approximate intermediate measure updates, adaptive tempering with an appropriate combination of filters is desirable and is the main focus of this work. One method for achieving this is introduced in \Cref{sec:tempering}, and extensively numerically explored in \Cref{sec:numeric}. In the following section, we will provide brief introductions to representatives of important filtering classes such as Gaussian-approximative filters, transport-based particle filter and McKean-Vlasov filters that we will amongst others consider for our approximate intermediate updates. 
\section{Gaussian-approximative filters}
The classical Kalman filter provides the optimum filtering solution for Gaussian noise and linear maps $\Psi$ and $H$, e.g. \cite{Kalman1960,Jazwinski2007}. 
In general, if either $\Psi$ or $H$ are nonlinear or the initial distribution is not Gaussian, it is not necessarily possible to identify the (parametric) family of distributions associated with the posterior. Empirical estimates, linearisation, or a combination of both are usually used to enforce a Gaussian approximation of the associated densities. 
A summaries of the different variants of Kalman based filters, e.g. Ensemble Kalman Filter (EnKF), Ensemble Square Root Filter, are given in \cite{Reich2015,Law2015, Evensen2022}.
The EnKF \cite{Evensen1994,Katzfuss2016} and its variants use the Kalman filter update for the mean and covariance to the their empirical approximation provided by the ensemble. In real-world applications, a stable implementation with small numerical error propagation is desirable. The stochastic noise of the vanilla variante of the EnKF \cite{AnalysisSchemeintheEnsembleKalmanFilter} can cause numerical challenges and is only producing correct moments in the linear case in the ensemble limit. As an alternative the family of Ensemble Square Root Filter was developed. These filters use rank reduction methods, and a square root matrix to overcome the matrix oversize problem, c.f. \cite{Tippett2003} and yield correct moments in the Gaussian case. In this article, we consider a particular variant within the family ensemble square root filters as described in \cite{Reich2015}, yet other members can be used interchangeably.

\subsection{Ensemble Square Root Filter}\label{sec:ESRF}
The Ensemble Square Root Filter (ESRF) as presented in \cite{Reich2015} utilises Monte Carlo estimates to  compute the forecast mean and anomalies as follows:
\begin{align}
    \overline{z}^f &=  \frac{1}{N_{ens}} \sum_{i=1}^{N_{ens}} z^f_{i}, \\
    A^f &= (z^f_{i} -\overline{z}^f).
\end{align}
The rank reduction approach, which involves computing only the anomalies of the covariance matrix instead of the full matrix, is used to prevent oversized vectors and matrices in high-dimensional systems. 
The main idea of the ESRF is to define the analysis step so that the empirical covariance of the resulting particle ensemble exactly matches the Kalman mean and covariance in the linear case:
\begin{align}
    P^f - KHP^f = P^a &= \frac{1}{N_{ens}-1} A^f S S^T (A^f)^T. \label{eq:SR_Pa}
\end{align}
Rearranging \eqref{eq:SR_Pa} leads to the following representation of the square root matrix $S$.
\begin{equation}
    S = \left( I_{N_{ens}} + \frac{1}{N_{ens}-1} (HA^f)^T R^{-1}  (HA^f)  \right)^{-1/2} \label{eq:SRM}
\end{equation}
ESRF by construction avoids loss of positive definiteness of the error covariance matrices.
Note that not all matrix square roots bear the desired relationship to the forecast ensemble: the analysis ensemble mean may not be equal to the analysis state estimate and consequently there may be an accompanying shortfall in the spread of the analysis ensemble as expressed by the ensemble covariance matrix. This points to the need for a restricted version of the notion of an ESRF, which is called an unbiased ESRF, c.f. \cite{Livings2008}. There is a set of generic necessary and sufficient conditions for the ensemble to remain centred on the analysis state estimate (unbiased). 
Hence, the entries $t_{ij}$ of the symmetric matrix $S$ should satisfies the following conditions:
\begin{equation}
    \sum_{i=1}^{N_{ens}} s_{ij} =\sum_{j=1}^{N_{ens}} s_{ij} =1.
\end{equation} 
Next one calculates the weights, which can be interpreted as importance weights, whereby $e_i$ denotes the i$^{\text{th}}$ basis vector in $\R^{N_{ens}}$: 
\begin{align}
    \hat{w}_i &=\frac{1}{N_{ens}}- \frac{1}{N_{ens}-1}  e_i^T S^2 (HA^f)^T R^{-1} (H\overline{z}^f-y^{obs}), \label{eq:ESRF_weight}
\end{align} 
with filter coefficients
\begin{align}
    d_{ij} &= \hat{w}_i - \frac{1}{N_{ens}} + s_{ij}. \label{eq:ESRF_filtercoeff}
\end{align}
The weights $\hat{w}_i$ and the coefficient $d_{ij}$ satisfy
\begin{equation}
    \sum_{i=1}^{N_{ens}} \hat{w}_i =1 \ \text{ and } \ \sum_{i=1}^{N_{ens}} d_{ij} =1.
\end{equation}
Finally, we can update the ensemble members with above equations using a general linear transformation from forecast to analysis state as described in \cite{Reich2015}, which is given by
\begin{equation}
    z_{j}^{a} = \sum_{j=1}^{N_{ens}} z^f_{i} d_{ij}.
\end{equation}
The implementation scheme can be found in the Appendix \Cref{alg:ESRF}. According to \cite{Chustagulprom2016} the computational complexity scales like $\mathcal{O}(N_{ens}^3)$ in the ensemble size.

\section{Approximative consistent filters}\label{sec:PF}
In the context of nonlinear evolution models, a non-parametric empirical approach, such as a particle filter (PF), is considered state-of-the-art for state and parameter estimation tasks, c.f. \cite{Fearnhead2018}, if sufficiently many particles can be generated. A key characteristic of these empirical filters is their consistency: they converge towards the true posterior density in nonlinear and non-Gaussian scenarios without requiring significant assumptions or restrictions on the prior or observations. Generally, particles of the prior are updated using individual importance weights derived from the likelihood. In high-dimensional spaces, particles are less likely to reside in regions of significant probability mass relative to the likelihood. Consequently, after performing importance sampling with respect to the likelihood, the posterior weights exhibit significant fluctuations. Many particles end up with very low relative weights, leading to a phenomenon known as filter collapse after a few iterations. This occurs because these low-weight particles are effectively lost for subsequent forecasting steps, demonstrating the filter's vulnerability to the curse of dimensionality. We characterise filter collapse by the effective sample size falling below a given threshold, as discussed in \cite{Frei2013}. The desire to mitigate these effects has inspired the development of various families of consistent filters. While none of the proposed filters are immune to the curse of dimensionality, they still do not assume any specific parametric structure for the underlying densities. However, it is important to mention that some form of approximation may still be involved (details in the respective sections).

\subsection{Bootstrap Particle Filter}\label{sec:Bootstrap}
The vanilla Bootstrap Particle Filter belongs to the family of Sequential Monte Carlo methods and is a simulation-based methods for calculating approximations to posterior distributions, c.f. \cite{Liu2001,Chopin2020}. As mentioned above this is achieved via importance weights, which are defined by: 
\begin{align}
    \hat{w}_{i} &\propto \exp \left( - \frac{1}{2}(Hz^f_{i}- y_k^{obs})^\top  R^{-1} (Hz^f_{i}- y_k^{obs}) \right) \cdot w_0, \ \ \ w_0=\frac{1}{N_{ens}}, \label{eq:Boot_weight} \\
    w &= \frac{\hat{w}_{i}}{\sum_{i=1}^{N_{ens}} \hat{w}_{i}}.
\end{align}
To address the issue of weight degeneracy, a resampling procedure is employed. This method of sequential importance resampling selectively removes particles with small weights and duplicates particles with large weights, effectively mitigating the degeneracy problem. Consequently, a new set of ensembles is resampled based on the calculated weights, resulting in an updated ensemble. However, this procedure can produce identical or near identical analysis ensemble members in practice. As a result, the filter can collapse. To overcome this issue, tuning tools are required c.f. \cite{Akyildiz2020, Cotter2020}. 
We present some of the variants thereof in \Cref{sec:add-on}. A pseudo-code implementation of the method is provided in the Appendix \Cref{alg:Bootstrap}.

\subsection{Transport based Particle Filter}
Transport-based filters approximate a coupling between the prior and posterior distributions. Typically, the goal is to find couplings that optimise an underlying optimal transport problem with respect to a given distance metric, such as the Monge-Kantorovich problem with the Wasserstein distance. It can be demonstrated that these optimal couplings are induced by a transport map. Several works \cite{ElMoselhy2012,Taghvaei2016,spantini2022, Cheng2013,Papamakarios2021} suggest approximative discrete optimal transport maps as surrogates of the true map. The advantage of using a map rather than a resampling step is that particles are moved rather than replicated, making them less dependent on the existing particles. However, it is important to note that the map is approximated based on the existing particles, which can potentially negatively affect the quality of the estimated map.  

\subsubsection{Ensemble Transform Particle Filter}\label{sec:ETPF}
The Ensemble Transform Particle Filter (ETPF) relies on a linear programming construction to approximate the optimal transport map, c.f. \cite{Reich2013}. As indicated above, the resulting linear transformation replaces the resampling step typically used in a standard Sequential Importance Resampling (SIR) filter\cite{Reich2015}. The map is approximated as follows
\begin{equation}
\label{eq:optimal_transport}
    T^\ast = \mathop{\mathrm{argmin}}
    \limits_{\substack{ T\in \mathcal{T} }} 
    \sum_{i,j=1}^{N_{ens}} t_{ij} d(z^f_{i}, z^f_{j})^2, 
\end{equation}
where  $\mathcal{T} = \{T\in \mathbb{R}^{N_{ens} \times N_{ens}} \colon
    \sum_{i=1}^{N_{ens}} t_{ij}=\frac{1}{N_{ens}} \land 
    \sum_{j=1}^{N_{ens}} t_{ij}=\hat{w}_{i} \}$ and $d$ is a suiteable metric on the state space i.e. on $\mathbb{R}^{N_z}$.
Note that the weights $\hat{w}_{i}$ are given like \eqref{eq:Boot_weight}. The implementation scheme can be found in the Appendix \Cref{alg:ETPF}. 
The ETPF is computationally expensive since one needs to solve the linear transport problem \eqref{eq:optimal_transport} in every assimilation step, c.f. \cite{Reich2015}.
Using the same algorithm to solve the system as in \cite{Chustagulprom2016}, i.e. FastEMD, the computational complexity scales like $\mathcal{O}(N_{ens}^3 \log(N_{ens}))$ in the ensemble size. As an alternative to the exact solution of the linear transport problem, one can compute the corresponding Sinkhorn approximation, which is obtained by regularising the underlying optimisation problem, c.f. \cite{Cuturi2013}. As the spread is severely underestimated, one can also compute a second-order correction extension of the ETPF \cite{Acevedo2017}.
\subsection{McKean-Vlasov Filters}
The general idea of moving particles instead of resampling them is also applied here. In particular, a modified evolution equation, i.e. in the continuous setting, is developed, incorporating the innovation term, i.e. describing the discrepancy between observation and projected state, weighted by a generalised gain function. This evolution equation directly generates samples from the posterior. To achieve this, the Fokker-Planck equation of the modified evolution equation is equated with the Kushner-Stratonovich equation, which describes the evolution of the posterior density in the continuous filtering problem. As a result, the generalised gain must be a solution of a weighted Poisson equation. Using this general approach, different families of filters can be derived. For a derivation of the general class of McKean-Vlasov Filters, see \cite{Pathiraja2021}. Our focus is on a variant of these filters, known as the Feedback Particle Filter (FPF) \cite{Yang2011,Taghvaei2018}.

\subsubsection{Feedback Particle Filter}\label{sec:FPF}
As previously stated, McKean-Vlasov Filters seek to derive a generalised gain that is utilised to update the prior sample with a structure similar to Ensemble Kalman filters, resulting in samples from the posterior \cite{Yang2011}:
\begin{equation}
    z_i^a = z_i^f - K_i^{(\epsilon,N_{ens})} \left( \frac{1}{2} (H(z_i^f)) + H(\overline{z}^f)) -y_{obs}  \right). \label{eq:FPF_analysis}
\end{equation}
To compute the gain $K$, a weighted Poisson equation 
\begin{align}
\begin{split}
\label{eq:Poisson}
       \nabla  \cdot (p(z,k) \nabla \phi(z,k)) &= - R^{-1}(H(z)- \hat{H})p(z,k), \\
    \int \phi(z,k) p(z,k) dz &=0, 
\end{split}
\end{align}
with $p$ being the conditional distribution of the state given the observation and $k \geq 0$, must be solved at every time step. Further, $R$ is the covariance matrix from the data model \eqref{eq:datamodel}.
The gain function is given for every time step $k$ by 
\begin{align}
    K(z,k) &= \nabla \phi(z,k), \label{eq:gain}
\end{align}
whereby $\phi$ is the solution of the Poisson equation. \cite{Yang2013} pointed out that $K$, in principle is, not uniquely determined. Yet \cite{Taghvaei2016} presented results on existence, uniqueness and regularity of \eqref{eq:Poisson} based on a weak formulation of the Poisson equation. Hence, the choice of the gain function in the multivariable case requires careful consideration of the uniqueness of the solutions of \eqref{eq:Poisson}. A significant challenge in this context is that \eqref{eq:Poisson} contains the posterior density, which is only available in the form of a sample approximation. Unfortunately, the estimation of the gain is not straightforward, since there is no closed-form solution of \eqref{eq:Poisson}, and the density $p$ is not explicitly known. Due to the limited amount of solution strategies available for high dimensional elliptic problems, we focus on a diffusion map \cite{Coifman2005,Coifman2006} approximation, theoretically investigated in \cite{Pathiraja2021a}. The underlying idea is to use diffusion map coordinates, which are obtained by uncovering the underlying geometry spanned by the particles. This approach transforms the problem into a lower-dimensional space, making it easier to solve. By leveraging the geometric structure of the data, the diffusion map approximation provides an efficient way to handle the complexities of the weighted Poisson equation within the McKean-Vlasov Filter framework. The first step is to construct a Markov matrix $T^{(\epsilon)}$ for $\epsilon>0$ defined on the state space approximating the semi-group, c.f. \cite{Grigor2006}, generated by $\frac{1}{p} \nabla \cdot ( p \nabla)$. This semi-group then provides the solution of the corresponding heat equation, which is given by
\begin{equation}\label{equ:heat}
    \frac{d}{dt}u  = \frac{1}{p} \nabla \cdot ( p \nabla) u.
\end{equation} 
Next we consider solutions to the initial value problem posed by \eqref{equ:heat} and $u(0)=\phi$ and denote them by $\phi^\epsilon = T^{(\epsilon)} \phi$.  
This, in turn, is used to formulate the following fixed-point problem derived (formally) from the fundamental theorem of calculus and \eqref{eq:Poisson} 
\begin{equation}
    \phi^{(\epsilon)} = T^{(\epsilon)} \phi^{(\epsilon)} + R^{-1} \int_0^{(\epsilon)} T^{(s)} (H-\hat{H}) ds, \label{eq:FPF_fixpoint}
\end{equation}
with
\begin{align}
    T^{(\epsilon)} \phi(z_i) &:= \frac{\int k^{(\epsilon)}(z_i,z_j) \phi(z_j) d\mu(z_j)}{\int k^{(\epsilon)}(z_i,z_j) d\mu(z_j)}, \\
    k^{(\epsilon)}(z_i,z_j) &:= \frac{g^{(\epsilon)}(z_i,z_j)}{\sqrt{\int g^{(\epsilon)}(z_i,z_j) d\mu(z_j)}\sqrt{\int g^{(\epsilon)}(z_i,z_j) d\mu(z_i)}}, \\
    g^{(\epsilon)}(z_i,z_j) &= \exp \left(  -\frac{\vert z_i-z_j\vert^2}{4 \epsilon} \right), \label{eq:FPF_g}
\end{align}
whereby $z_i$ and $z_j$ are ensembles of the state vector, $\mu$ is an absolute continuous probability measure with respect to the posterior $p$, and $g$ is the Gaussian kernel. 

For a numerical implementation one uses $\epsilon (H-\hat{H}) \sim \int_0^{\epsilon} T^{(s)} (H-\hat{H}) ds$.
The particle approximation is obtained by Monte Carlo integration and empirical approximation of $T^{(\epsilon)} \phi$. To the knowledge of the authors, there is a gap in the literature when implementing the FPF in time-discrete scenarios. For the time continuous case, $z_i$ are already samples consistent with the filtering distribution, i.e. the posterior. However, in the time discrete case, only the prior samples $z_i^f$ are available. In general, they are not consistent with $\mu$, for every $\Delta t_{obs}>0$. 
Therefore, we propose to use an empirical approximation, to our understanding different from \cite{Yang2011,Yang2013,Taghvaei2018,Taghvaei2020}.
To this end, we apply Monte Carlo integration, importance sampling, and \eqref{eq:filtering_problem}
\begin{equation}
\begin{aligned}
    T^{(\epsilon)} \phi(z_i) & =  \frac{\int k^{(\epsilon)}(z_i,z_j)\phi(z_j) d\mu(z_j)}{\int k^{(\epsilon)}(z_i,z_j) d\mu(z_j) }\\
    &= \frac{\sum_j k^{(\epsilon)}(z_i,z_j) \Phi^{(\epsilon,N_{ens})}_j e^{- \frac{1}{2}(H(z_j)-y_{obs})^T R^{-1}(H(z_j)-y_{obs})}}{ \sum_j k^{(\epsilon)}(z_i,z_j) e^{- \frac{1}{2}(H(z_j)-y_{obs})^T R^{-1}(H(z_j)-y_{obs})} }
\end{aligned}
\end{equation}
and subsequently state the following empirical approximation of the Markov transition matrix
\begin{equation} \label{eq:FPF_T_MC}
 T^{(\epsilon,N_{ens})}_{i,j}=  \frac{k^{(\epsilon)}(z_i,z_j) e^{- \frac{1}{2}(H(z_j)-y_{obs})^T R^{-1}(H(z_j)-y_{obs})}}{ \sum_j k^{(\epsilon)}(z_i,z_j) e^{- \frac{1}{2}(H(z_j)-y_{obs})^T R^{-1}(H(z_j)-y_{obs})} }.
\end{equation}
Note that $\Phi^{(\epsilon,N_{ens})}_j = \phi(z_j)$ for $j=1,...,N_{ens}$.
A pseudo implementation scheme is given in Appendix \Cref{alg:Feedback}. 
It is important to state, that despite the consistency theory, provided in \cite{Taghvaei2020}, the numerical resolution of \eqref{eq:FPF_fixpoint} is unstable for fixed particle number $N_{ens}$ and $\epsilon \rightarrow 0$. 
This is not a contradiction to other well-posedness statements in \cite{Taghvaei2020}, as the corresponding limits in $\epsilon$ and $N_{ens}$ do not  commute.
When considering \eqref{eq:FPF_fixpoint} and \eqref{eq:FPF_g}, it comes at no surprise that depending on the dimension of the state space, we must scale the bandwith of the kernel accordingly, otherwise $T^{(\epsilon)}$ coincides with the identity in every reasonable finite floating point precision and \eqref{eq:FPF_fixpoint} becomes ill-posed.
For the actual implementation in our time discrete scenario, we adapt this estimation
To avoid another tuning parameter, we scale the bandwidth/time parameter with the maximum distance between particles. This obviously leads to higher approximation errors of the integral term, but considering the arguably low number of particles we do not expect this to be a relevant source of approximation errors. However, one could improve the latter by using better quadrature formulas to approximate the integral.

\section{Adaptive tempering}\label{sec:tempering}
Here, we propose an adaptive approximative tempering method, which is a variant of the traditional tempering introduced in \Cref{sec:math}. Updates of the intermediate measures \eqref{eq:intermediate_measures} are performed using members of different families of filters, including approximative Gaussian filters. This approach introduces an error, as the true posterior cannot be fully preserved within the tempering procedure. However, the advantage of this approximative ansatz lies in its enhanced robustness in scenarios where individual representatives of approximative consistent filter families encounter difficulties. For our tempering schedule, we restrict ourselves to \( T = 2 \), which corresponds to \(\tau_0 = 0\), \(\tau_1 = \alpha\), and \(\tau_2 = 1\). While this method can be extended to include more tempering steps (either through traditional consistent tempering or approximative updates), the increased number of possible combinations significantly raises the tuning overhead, which is already non-trivial in the simplified case considered here. Using the notation in \eqref{eq:intermediate_measures} for our simplified schedule yields the following tempered likelihood: 
\begin{align}
\begin{split} \label{eq:split_L}
    p({\color{black}y_{k}} \vert {\color{black}z_{k}}) &\propto \exp \left( \frac{\alpha}{2}(Hz^f_{i,k}- y_k^{obs})^\top  R^{-1} (Hz^f_{i,k}- y_k^{obs}) \right)  \\
    &\times \exp \left( \frac{1-\alpha}{2}(Hz^f_{i,k}- y_k^{obs})^\top  R^{-1} (Hz^f_{i,k}- y_k^{obs}) \right) .
\end{split}    
\end{align}
This or similar splitting approaches have been explored in the literature for specific filters \cite{Chustagulprom2016, Frei2013}. This strategy is often referred to as a hybrid approach, aiming to combine filters with complementary properties. Typically, this involves selecting one filter from the family of Gaussian approximative filters and one from the family of (approximative) consistent filters. We introduce key representatives from the important classes of consistent filters and focus on the family of Ensemble Square Root Filters within the Gaussian approximative filter family.
While it is possible to consider different filters from these two classes, the key lies in the ability to adjust the likelihood in the filter under consideration. Since the likelihood is incorporated differently in the filter family examined in this paper, we will provide a brief overview of how the split affects the updates.
\subsection{Tempered ESRF}
For the ESRF, the square root matrix \eqref{eq:SRM} and the weights \eqref{eq:ESRF_weight} change into:
\begin{align}
\label{eq:ESRF_temp_S}
S(\alpha) &= \left( I_{N_{ens}} + \frac{1-\alpha}{N_{ens}-1} (HA^f)^T R^{-1}  (HA^f)  \right)^{-1/2} \\
\hat{w}_i(\alpha) &=\frac{1}{N_{ens}}- \frac{1-\alpha}{N_{ens}-1}  e_i^T S^2 (HA^f)^T R^{-1} (H\overline{z}^f-y^{obs}). \label{eq:ESRF_temp_w}
\end{align}
\subsection{Tempered Bootstrap and ETPF}
For the Bootstrap and ETPF the weights in \eqref{eq:Boot_weight} change into:
\begin{align}
    \hat{w}_{i,k}(\alpha) &\propto \exp \left( - \frac{\alpha}{2}(Hz^f_{i,k}- y_k^{obs})^\top  R^{-1} (Hz^f_{i,k}- y_k^{obs}) \right) \cdot w_0. \label{eq:PF_temp_w}
\end{align}
Additionally, the transport map \eqref{eq:optimal_transport} from the ETPF changes as follows:
\begin{equation}
    T^\ast(\alpha) = \mathop{\mathrm{argmin}}
    \limits_{\substack{ T\in \mathcal{T} }} 
    \sum_{i,j=1}^{N_{ens}} t_{ij} d(z^f_{i}, z^f_{j})^2, 
\end{equation}
with the constraints
\begin{equation}
 \sum_{i=1}^{N_{ens}} t_{ij}=\frac{1}{N_{ens}}, \ \ \  
    \sum_{j=1}^{N_{ens}} t_{ij}=\hat{w}_{i}(\alpha).
\end{equation}
\subsection{Tempered FPF}
For the FPF we insert the bridging parameter in \eqref{eq:FPF_T_MC} as follows:
\begin{equation}
 T^{(\epsilon,N_{ens})}_{i,j}=  \frac{k^{(\epsilon)}(z_i,z_j) e^{- \frac{\alpha}{2}(H(z_j)-y_{obs})^T R^{-1}(H(z_j)-y_{obs})}}{ \sum_j k^{(\epsilon)}(z_i,z_j) e^{- \frac{\alpha}{2}(H(z_j)-y_{obs})^T R^{-1}(H(z_j)-y_{obs})} }.
\end{equation}
\subsection{Tuning algorithm }
At this point, it remains unclear in which filtering cycles tempering is beneficial and how to choose the bridging parameter \(\alpha\). One brute-force approach to selecting the value for \(\alpha\) is to evaluate the bridging parameter in every assimilation cycle, as initially discussed in \cite{Frei2013}. Following this approach, \cite{Nerger2021} provides a detailed study on the choice of the bridging parameter for the tempered combination of the local ensemble transform Kalman filter with the nonlinear ensemble transform filter. Traditionally, the effective sample size is considered state-of-the-art for adjusting the tempering schedule. While these techniques can easily be combined with our adaptive tempering approach, we limit ourselves to fixed values of \(\alpha \in \{0.2, 1\}\). This choice simplifies the evaluation of the different combinations already being considered and aligns with the intuition that small updates with beneficial properties are sufficient to adjust the main filter being used. Additionally, this choice is consistent with the values found to be optimal by numerical experiments for a complementary filter combination considered in \cite{Chustagulprom2016}. In previous works on the type of tempered filters considered in here, tempering was fixed to occur in every data assimilation cycle.
Our main focus is to adaptively decide in each filtering step whether a tempering step is required. Intuitively, the effective sample size offers a measure of the number of particles actively contributing to the update procedure. This approach is not only intuitively justified but also provides valuable statistical insights for tempering and the decision process of whether to resample. We aim to establish a even more incisive and informative criterion for determining whether an approximative tempering step is necessary. The novel criterion we propose relies on general descriptive statistics of the particles mapped to the observational space. To benchmark against the state-of-the-art approach, we compare our numerical results with the fixed tempering method employed in earlier works on hybrid filtering (for details see \Cref{sec:numeric}).

\subsubsection{Tempering criterion I}\label{sec:ESS}
In a general particle filter setting the effective sample size (ESS) is used as an effective criterion to resample and is given as 
\begin{equation}
    ESS = \frac{1}{\sum_{i=1}^{N_{ens}} w_i^2}.
\end{equation}
This approach simplifies the comparison of filters that require an update of the particles at each step. Rather than using the ESS for resampling, we use it to decide if tempering, which entails two complementary filtering updates, should be performed. Specifically, the default is to perform updates only with one filter, unless the ESS drops below a certain threshold \(\theta\).
\begin{equation}
     ESS < \theta N_{ens}.
\end{equation}
The threshold is set at \(\theta = 0.5\) times the ensemble size in this case, which is a commonly suggested value in the resampling context. Naturally, the hyperparameter \(\theta\) can also be adjusted to more extreme values and this is certainly a topic of interest for future research.
\subsubsection{Tempering criterion II}\label{sec:IQR}
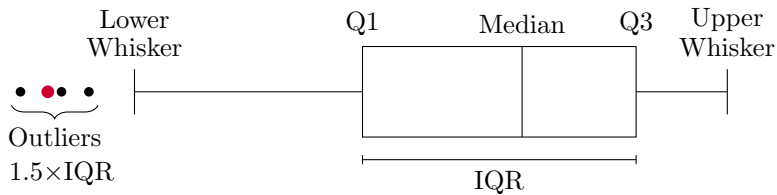
\begin{figure}[ht]
    \centering
\begin{tikzpicture}[scale=0.6]
    \draw (0,0) rectangle (6,2);
    \draw (3.5,0) -- (3.5,2);

    \draw (0,-0.5) -- (6,-0.5);
    \draw (0,-0.4) -- (0,-0.6);
    \draw (6,-0.4) -- (6,-0.6);

    \node at (3.5,2.5) {Median};
    \node at (0,2.5) {Q1};
    \node at (6,2.5) {Q3};
    \node at (-5,2.6) {Lower};
    \node at (-5,2.0) {Whisker};
    \node at (8,2.6) {Upper};
    \node at (8,2.0) {Whisker};
    \node at (3,-1) {IQR};
    \node at (-6.6,-0.8) {1.5$\times$IQR};

    \draw (0,1) -- (-5,1);    
    \draw (-5,1.5) -- (-5,0.5);

    \draw (6,1) -- (8,1);
    \draw (8,1.5) -- (8,0.5);

    \fill (-6,1) circle (3pt);
    \fill (-6.6,1) circle (3pt);
    \fill (-7.5,1) circle (3pt);
    \fill[cherryred] (-6.9,1) circle (4pt);
    
    \draw[decorate, decoration={brace, amplitude=5pt, mirror}] (-7.7,0.7) -- (-5.8,0.7) node[midway, below=5pt] {Outliers};
    
\end{tikzpicture}
    \caption{Illustration of a typical box plot. The \textcolor{cherryred}{red} dot denotes the observation. }
    \label{fig:Boxplot_tempering}
\end{figure}

Since the decision-making process should be based on the discrepancy between observed particle positions and the observational data, we propose using a criterion grounded in descriptive statistics for this property. Specifically, we apply an individual filter of our choice if the current observation falls within the interquartile range (IQR) of the observed particles. If the observation is outside the IQR, we employ a tempered combination approach: one update is performed with the currently considered filter, and the other update is performed with a filter from a different family.
The main aim is to use consistent filters to enhance the uncertainty quantification of the estimate. These filters typically require periodic adjustment through robust alternative filters. However, we also consider Gaussian approximative filters as primary filters combined with various other filters to investigate the associated benefits.
 An illustration of a box plot is given in \Cref{fig:Boxplot_tempering}, whereby the observation is outside the IQR, which illustrates the case when tempering is applied. Box plots and the underlying descriptive statistics offer an effective graphical representation of data concentration, as illustrated in \Cref{fig:Boxplot_tempering}. They also accurately depict the extent to which extreme values deviate from the central tendency of the data. A box plot represents five key statistics: the minimum (lower whisker), the first quartile (Q1), the median (Q2), the third quartile (Q3), and the maximum (upper whisker). Further, outliers are marked as points as depicted in \Cref{fig:Boxplot_tempering}. The position of the quantiles are determined via:
\begin{equation}
    Q_i = \frac{i \cdot (N_{ens}+1)}{4}  \ \text{ for } i = 1,2,3. 
\end{equation}
An outlier is a data point that differs numerically from the rest of the data. They are determined by a  quantile factor $1.5$ times the interquartile range (IQR), which is given as
\begin{equation}
    IQR(\{\xi_i\}_{i \in\{ 1, \dots ,N_{ens}\}}) = [\xi_{\sigma(Q_3)}, \xi_{\sigma(Q_1)}].
\end{equation}
where \(\xi_i \in \mathbb{R}\) is an ensemble of real numbers and $\sigma$ is the permutation such that $\xi_{\sigma(i)}$ is monotonously ordered for $i\in \{1,\dots, N_{ens}\}$.
The IQR interval relates to fifty per cent of the score. The usual value of the quartile factor, associated with detecting minor outliers of the sample set, is $1.5$. Hence, we differ the length of the IQR with the quantile factor $\lbrace 0,1.5\rbrace$. 
Note that the lower and upper whisker excludes the outliers.
Twenty-five per cent of scores fall below the $Q1$. Contrary, Twenty-five percent of data are above $Q3$. Meaning that seventy-five per cent of the scores fall below the third quartile.
\Cref{fig:Boxplot_tempering} depicts a line that divides the box into two portions to indicate the median, which represents the midpoint of the data. Half of the scores exceed or equal this value, while the other half fall below it. The box plot shape indicates whether a statistical dataset is normally distributed or skewed. 
When the median is in the centre of the box and the whiskers are roughly equal on both sides, the distribution is symmetric.
When the median approaches $Q1$ and the whisker on the lower end of the box is shorter, the distribution is positively skewed, i.e. skewed right.
When the median is closer to $Q3$ and the whisker is shorter on the top end of the box, the distribution is negatively skewed, i.e. left skewed.
Furthermore, the box plot is useful as it shows the dispersion of a data set. The data is more evenly distributed as the IQR length increases. The smaller the IQR, the lower the spread of the data. One can determine how similar other data values are.
To apply the criterion to an ensemble of particles in $z_i \in \mathbb{R}^{N_z}$ we first observe the states and subsequently project onto each individual component of the observation i.e. 
we compute 
\begin{equation}
IQR(\{e_j^T Hz_i\}_{i\in \{1,\dots ,N_{ens}\}})
\end{equation}
for each component $j\in \{1,\dots, N_{y}\}$ and with $e_j$ being the canonical basis vector in $\mathbb{R}^{N_{y}}$.
Subsequently, the criterion then reads as in \Cref{alg:iqr}.
\begin{algorithm}[H]
\caption{IQR criterion for ensembles in $\mathbb{R}^{N_{ens} \times N_y}$}
\label{alg:iqr}
    \begin{algorithmic}
        \If{ $e_j^T y_{obs} \in IQR(\{e_j^T Hz_i\}_{i\in \{1,\dots, N_{ens}\}}) \quad \forall j \in \{1, \dots, N_{y}\} $} 
        \vspace{0.25\baselineskip}
         \State $\alpha=1$ 
        \Else 
        \State $\alpha=0.2$
        \EndIf
    \end{algorithmic} 
\end{algorithm}

\section{Numerical Examples}\label{sec:numeric}
We demonstrate the numerical behaviour of the proposed combination of the aforementioned filters and the related complementary filtering techniques. 
Our choice of experiments includes benchmarks based on the classical models introduced in \cite{Lorenz1963,Lorenz1996} and the reference setup provided in \cite{Chustagulprom2016}. Additionally, we investigate the behaviour of existing and the proposed tempering criteria using the Langevin-type stochastic differential equation originating from a double well potential. We consider this scenario as prototypical for many situations in high-dimensional state spaces, whenever the (particle) filter is far from the asymptotic regime i.e. the number of particles is small in comparison to the state space and the regularity of the filtering distribution.  
Finally, we consider the one-dimensional shallow water equations to provide a descriptive illustration and physical interpretation of the proposed tempering criteria using a spatio-temporal model. 
Additionally, to the baseline of each individual filter, we test three different tempering approaches, c.f. \Cref{sec:tempering}. In the first one, we use the tempered combination in each data assimilation step, denoted by ETPF-ESRF. The second one includes the ESS as criteria to temper, i.e. ESS-ETPF-ESRF \Cref{sec:ESS}. Meaning that the tempered version is only used in case the $ESS < 0.5 N_{ens}$. Last but not least, we test the IQR approach, denoted by IQR-ETPF-ESRF, c.f. \Cref{sec:IQR}. Meaning that we only use the tempered filter combination in case the observation is not within the IQR. Before providing the first example, we would like to point out the crucial additional strategies used to ensure the filters are both robust and computationally practical.

\subsection{Add-on filter techniques}\label{sec:add-on}

 Due to computational constraints that limit ensemble size, a critical consideration is mitigating the effects of small ensembles. Specifically, filters are susceptible to sampling errors, which diminish accuracy and lead to underestimation of error covariances. This can result in filter divergence, where the filter fails to track the reference solution effectively. Techniques such as localization or inflation/rejuvenation can address these issues.

\subsubsection{Inflation}\label{sec:inflation}
Inflation, for instance, increases forecast uncertainty to enhance the robustness of ensemble filter implementations against the effects of finite ensemble sizes. 
Given a forecast ensemble $z_{i,k}^f$, we first compute its empirical mean $\overline{z_k}^f$ and the ensemble anomalies $A_k^f$. An inflation ensemble is defined according to \cite{Reich2015} by
\begin{equation}
    z^f = \overline{z}^f + \gamma A^f
\end{equation}
whereby $\gamma>1$ is the inflation factor. Inflation is used to avoid filter divergence for small ensemble sizes. 
However, inflation does not directly account for spatial structures.
\cite{Cotter2020} presented two further tuning tools related to inflation, namely jittering and nudging. Jittering improves the diversity of the ensemble by computing new ensembles, and nudging corrects the solution term of the evolution model to keep ensembles closer to the true state. An optimal nudging procedure in particle filters is given in \cite{Lingala2014}.

\subsubsection{Rejuvenation} \label{sec:rejunvenation}
A general sequential importance resampling filter can produce identical or near identical ensemble members.
This issue can be overcome according to \cite{Chustagulprom2016} by using particle rejuvenation in the following way:
\begin{equation}
    z_j^a = z_j^a +\sum_{i=1}^{N_{ens}} (z_i^f -\overline{z}^f ) \frac{\tau \xi_{ij}}{\sqrt{N_{ens}-1}},
\end{equation}
\noindent 
where $\tau>0$ denotes the rejuvenation parameter and $\xi_{ij}$ are iid Gaussian random variables with mean zero and standard deviation one. Rejuvenation is employed to prevent the occurrence of identical or nearly identical ensemble members, while also increasing the ensemble's spread.

\subsubsection{Localisation}\label{sec:localisation}
\begin{figure}[ht]
\centering
\begin{tikzpicture}[scale=0.66]

        \draw[step=1cm,gray,very thin] (0,0) grid (6,5);

        \foreach \x in {0,1,2,...,6}
            \foreach \y in {0,1,2,...,5}
                \fill[darkblue]  (\x,\y) circle (3pt);
                
         \draw[fill=darkgreen](2,2)circle(3pt);   
         
         \draw[darkgreen] (2,2) circle (1.25);

         \draw[fill=darkorange] (0.4,0.4) rectangle (0.8,0.8);
          \draw[fill=cherryred] (1.4,2.4) rectangle (1.8,2.8);
          
           \draw[fill=darkorange] (2.4,3.4) rectangle (2.8,3.8);

            \draw[fill=darkorange] (5.4,4.4) rectangle (5.8,4.8);

            \draw[fill=darkorange] (4.4,1.4) rectangle (4.8,1.8);

\end{tikzpicture}
\caption{ Illustration of localization, \textcolor{darkblue}{Blue} are the grid points, \textcolor{darkgreen}{green} the currently updated point, \textcolor{darkorange}{orange} the available observation and \textcolor{cherryred}{red} the selected obs}
\label{fig:Localisation}
\end{figure}
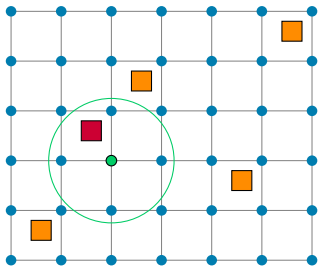

Difficulties arising when filtering spatio-temporal processes are twofold. Firstly the ensemble size is much smaller than the state space and secondly, spurious correlations are introduced as not all observations are physically relevant for each point in the domain. Both issues can at least partly mitigated by localisation, a technique, which determines the analysis state by conducting "local analyses" at each model grid point, considering only observations within a specified nearby region, i.e. \Cref{fig:Localisation}.
The seminal work \cite{Hunt2007} provides a standard implementation for this technique applied to Ensemble Kalman filter and is commonly used in numerical weather prediction and climate analysis, as it realistically aligns with the underlying dynamics. The general idea can be implemented in multiple ways. A prevalent technique involves modifying the ensemble's covariance matrix to mitigate spurious long-range correlations. Known as B-localization, this method is widely acknowledged yet challenging to apply uniformly across various filtering techniques, c.f. \cite{Reich2015}. Conversely, R-localization incorporates only observations within a predefined vicinity of the grid point under consideration. Similar to B-localization, the selection of the effective radius for local regions in R-localization is critical, as it must accurately reflect the spatial extent of significant dynamical correlations manifested by the ensemble. 
Considering the differences in how the R matrix and individual updates are adjusted across different filtering families, we will specify the localisation schemes employed for the families of filters addressed in this study.
According to \cite{Reich2015}, the localisation matrix $\Tilde{C_g} \in \mathbb{R}^{N_y \times N_y}$ with localisation radius $r_{loc}$ is given by:
\noindent
\begin{equation}
  (\Tilde{C_g})_{ll}= \rho \left( \frac{\vert x_g - \Tilde{x_l}\vert_L}{r_{loc}} \right), \label{eq:loc_matrix}
\end{equation}
\noindent
with the current grid point $x_g = g \Delta x$ and the observation grid point $\Tilde{x_l}=l \Delta x$. $\rho$ is known as filtering function and can have linear or higher polynomial order, e.g. \cite{Reich2015}. 
\noindent
\Cref{fig:Localisation} shows an example of the localisation technique described above. The green point is the current grid point. If an observation is within the radius of the current grid point, it is taken into account when calculating the localisation matrix. All other observations are not of interest. 

 \subsubsection{LESRF}
 The localised Ensemble Square Root Filter (LESRF) inserts the localisation matrix \eqref{eq:loc_matrix} directly into the square root matrix $S$, c.f. \cite{Reich2015}. In addition, each grid point can now be treated separately. Parallel processing can now speed up the algorithm.
Using inflation as well as localisation as a tuning tool, the algorithm of the ESRF, Appendix \Cref{alg:ESRF}, changes into:
 
\begin{algorithm}[H]
\small
    \caption{Localised Ensemble Square Root Filter (LESRF) (LETKF, \cite{Hunt2007}) see also the presentation in \cite{Reich2015}} \label{alg:LocESRF}
   \begin{algorithmic}[1]
        \vspace{0.25\baselineskip}
        \Require{ensemble forecast $z^f_i\in \mathbb{R}^{N_x}, \quad  \forall i\in \{1\dots N_{ens}\}$ and observation $y^{obs} \in \mathbb{R}^{N_y}, \gamma=1.05$}
        \vspace{0.5\baselineskip} 
        \State $\overline{z}^f =  \frac{1}{N_{ens}} \sum_{i=1}^{N_{ens}} z^f_{i}$
        \vspace{0.5\baselineskip}
        \State $A^f = (z^f_{i} -\overline{z}^f)$
        \vspace{0.5\baselineskip}
         \State inflation: $z_i^f = \overline{z}^f + \gamma A^f$
        \vspace{0.5\baselineskip}
        \State compute step 1 and 2 again
        \vspace{0.5\baselineskip}
        \ForAll{ $x_g \in \Omega_{N_z}$}
        \State $S[x_g] = \left( I_{N_{ens}} + \frac{1}{N_{ens}-1} (HA^f)^T \Tilde{C_g} R^{-1}  (HA^f)  \right)^{-1/2} $
        \vspace{0.5\baselineskip}
        \State $ \hat{w}_i[x_g] =\frac{1}{N_{ens}}- \frac{1}{N_{ens}-1}  e_i^T T^2[x_g] (HA^f)^T \Tilde{C_g} R^{-1} (H\overline{z}^f-y^{obs}) $
        \vspace{0.5\baselineskip}
        \State $ d_{ij}[x_g] = \hat{w}_i[x_g] - \frac{1}{N_{ens}} + t_{ij}[x_g] $
        \vspace{0.5\baselineskip}
        \State $z_{j}^{a}[x_g] =\sum_{i=1}^{N_{ens}}  z^f_{i}[x_g] d_{ij}[x_g] $
        \vspace{0.5\baselineskip}
        \EndFor\\
        \vspace{0.5\baselineskip}
        \Return $z^a_{j}$ for every $i\in \{1\dots N_{ens}\}$
   \end{algorithmic} 
\end{algorithm}

\subsubsection{LETPF}
For the localised Ensemble Transform Particle Filter (LETPF), a localised distance between two ensemble members at each grid point $x_g$ is defined as described in \cite{Chustagulprom2016}.
\begin{equation}
    \Tilde{d}[x_g](z_{i_1}^f ,z_{i_2}^f)^2 = \sum_{l=1}^{N_d} \rho \left( \frac{\vert x_g - \Tilde{x_l}\vert_L}{r_{loc}} \right) d_2\left( z_{i_1}^f(\Tilde{x_l}) , z_{i_2}^f(\Tilde{x_l}) \right)^2
\end{equation}
This distance is a weighted average of the differences between the two ensemble members over nearby gird points \cite{Reich2015}. Using rejuvenation as well as localisation as tuning tool, the algorithm of the ETPF, Appendix \Cref{alg:ETPF}, changes into:

\begin{algorithm}[]
\small
\caption{Localised Ensemble Transform Particle Filter (LETPF), see also the presentation in \cite{Reich2015}} \label{alg:LocETPF}
\begin{algorithmic}[1]
\vspace{0.25\baselineskip}
\Require{ensemble forecast $z^f_i\in \mathbb{R}^{N_x}, \quad  \forall i\in \{1\dots N_{ens}\}$ and observation $y^{obs} \in \mathbb{R}^{N_y}, \tau=0.2$}
\vspace{0.5\baselineskip}
\State $\overline{z_k}^f =  \frac{1}{N_{ens}} \sum_{i=1}^{N_{ens}} z^f_{i}
$
\vspace{0.5\baselineskip}
\State $A^f = (z^f_{i} -\overline{z}^f)$
\vspace{0.5\baselineskip}
\ForAll{ $x_g \in \Omega_{N_z}$}
\State  $ w_{i}[x_g] \propto \exp \left( - \frac{1}{2}(Hz^f_{i}- y^{obs})^\top  \Tilde{C_g} R^{-1} (Hz^f_{i}- y^{obs}) \right) $ 
\vspace{0.5\baselineskip}
\State $\hat{w}_{i}[x_g] = \frac{w_{i}[x_g]}{\sum_{i=1}^{N_{ens}} w_{i}[x_g]}$ 
\vspace{0.5\baselineskip}
\State $T^\ast[x_g] \gets$ solution of \eqref{eq:optimal_transport} using the localised metric i.e. \(d=\tilde{d}[x_g]\).  
\vspace{0.5\baselineskip}
 \State $rejuvenation= \frac{\tau }{N_{ens}-1} A^f \xi, \ \xi \sim \mathcal{N}(0,1)$
        \vspace{0.5\baselineskip}
\State $z^a_{i}[x_g] = \sum_{i=1}^{N_{ens}} z^f_{i}[x_g] N_{ens}t_{ij}^\ast[x_g] + rejuvenation$
\vspace{0.25\baselineskip}
\EndFor\\
\vspace{0.25\baselineskip}
\Return $z^a_{i}$ for every $i\in \{1\dots N_{ens}\}$
\end{algorithmic}
\end{algorithm} 

\subsection{Langevin Dynamics} \label{res:Langevin}
\begin{figure}[ht]
\centering 
\includegraphics[width=0.49\textwidth]{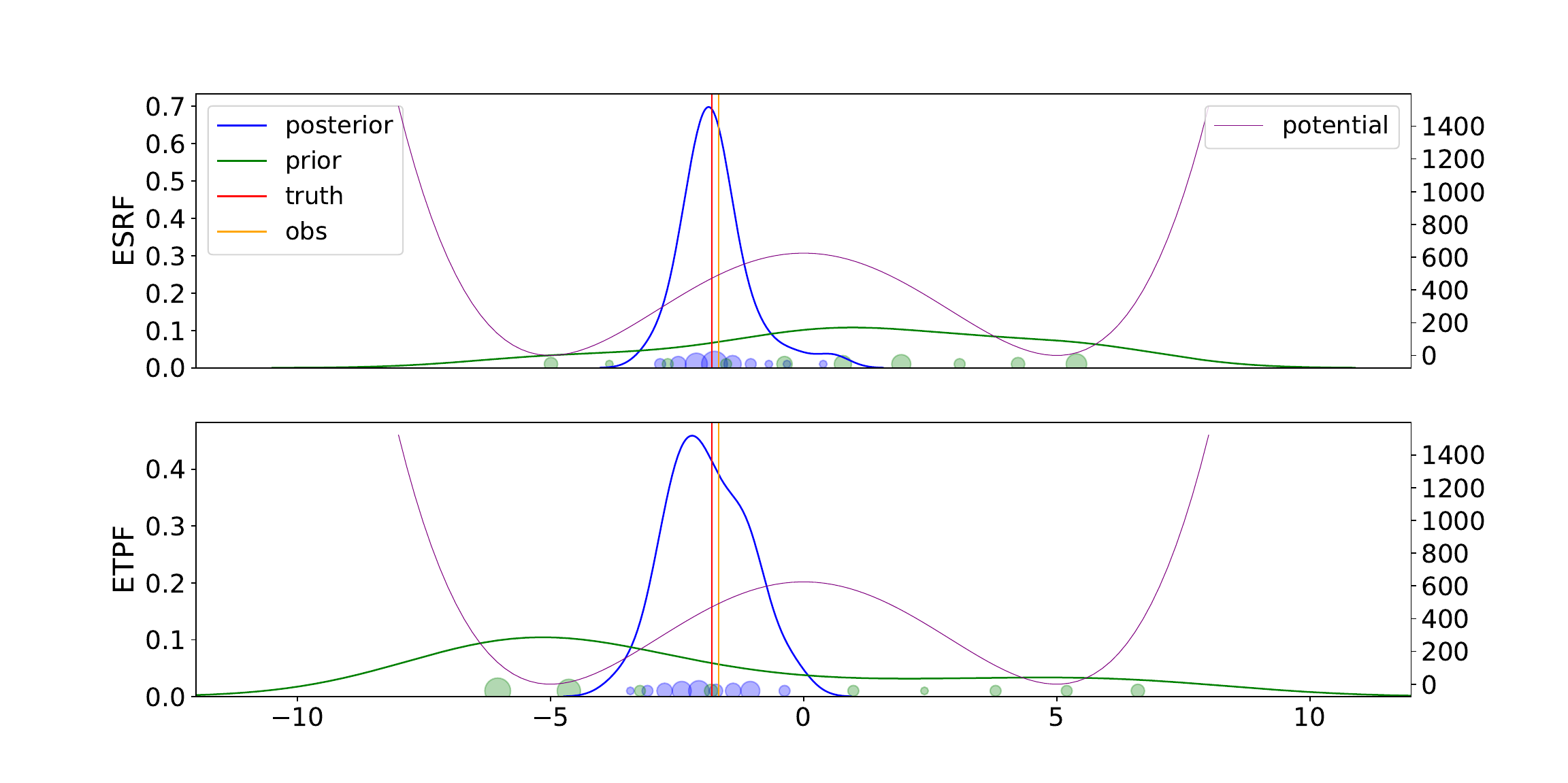}
\includegraphics[width=0.49\textwidth]{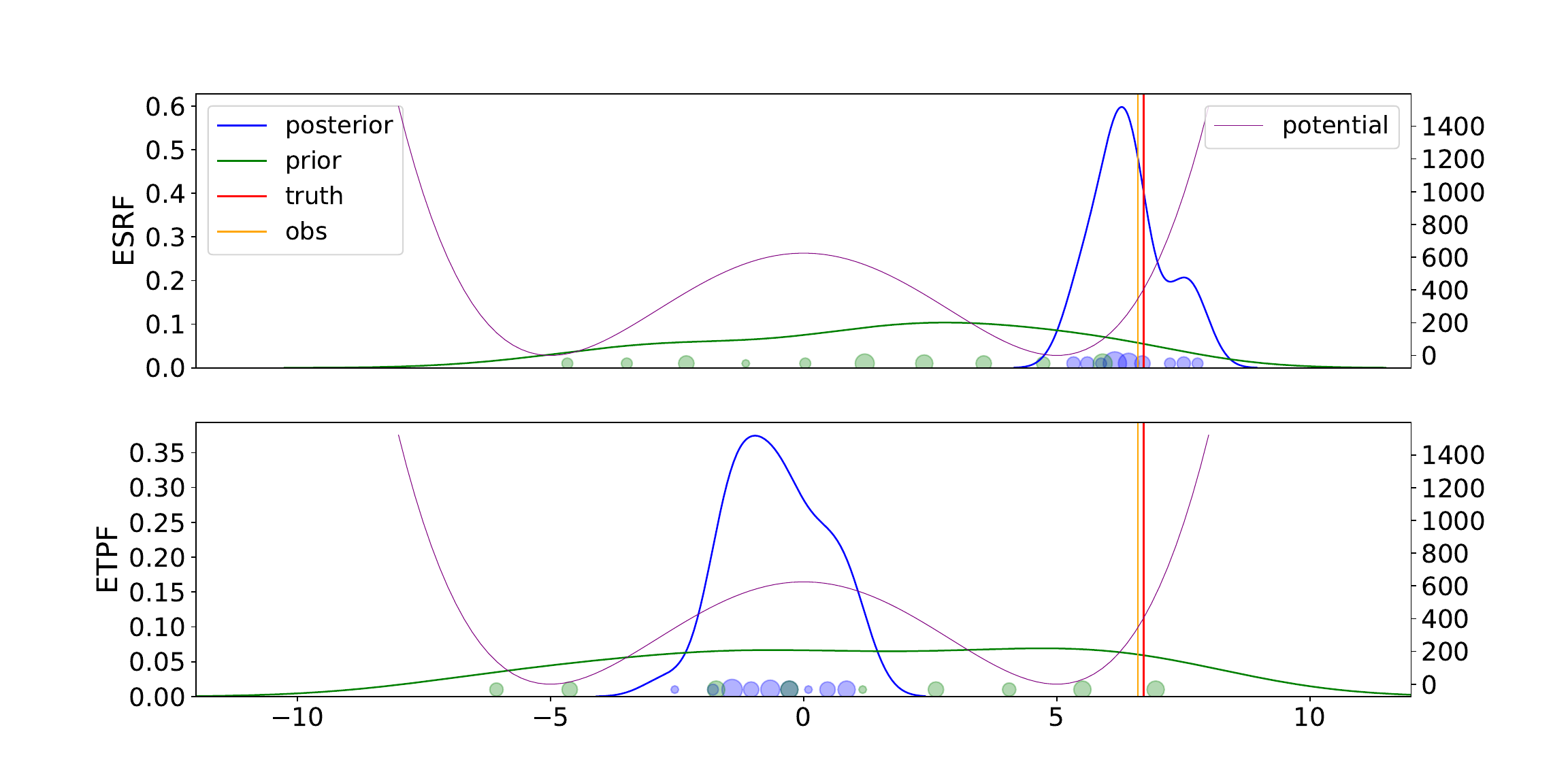}
\includegraphics[width=0.49\textwidth]{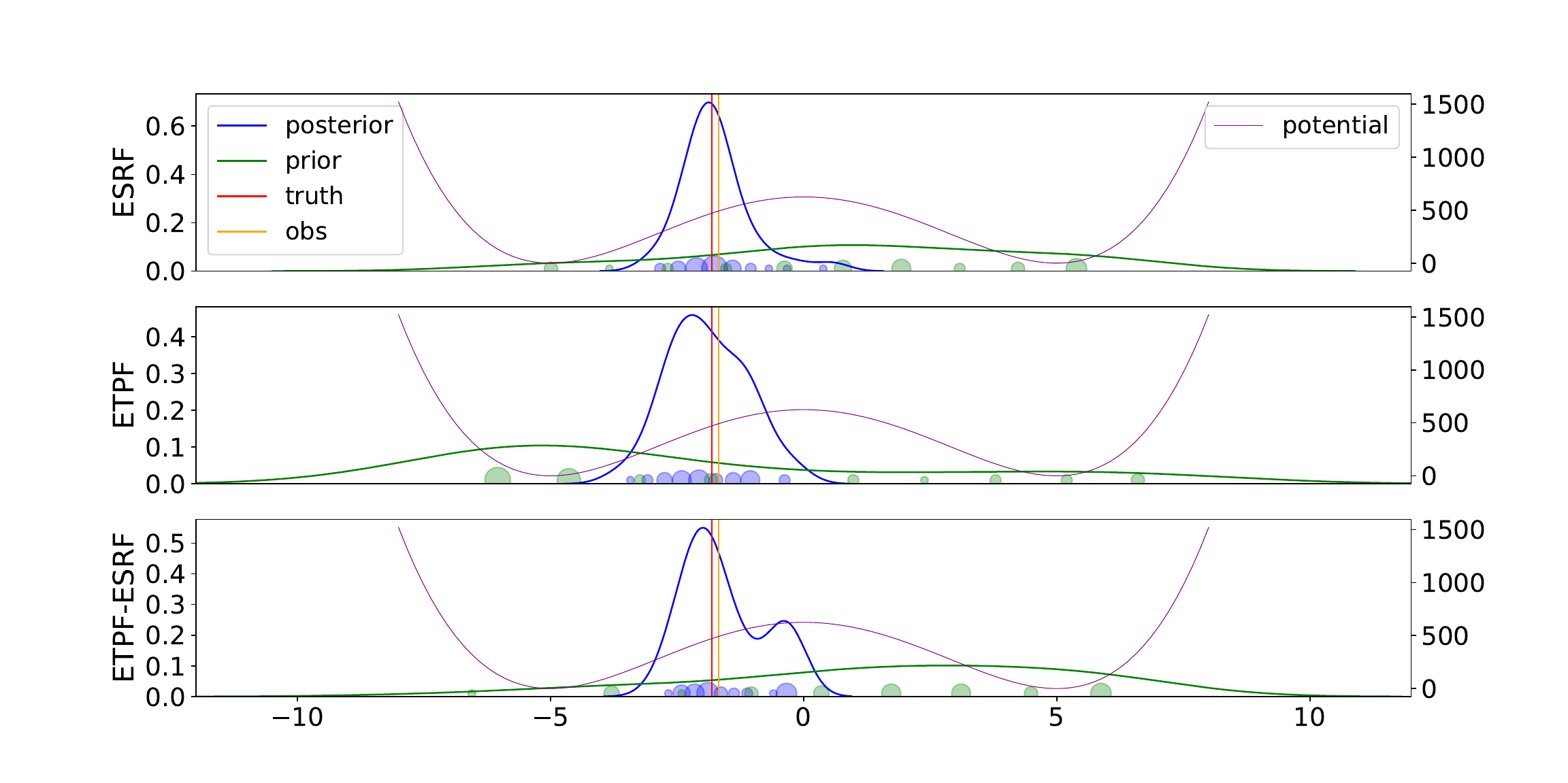}
\includegraphics[width=0.49\textwidth]{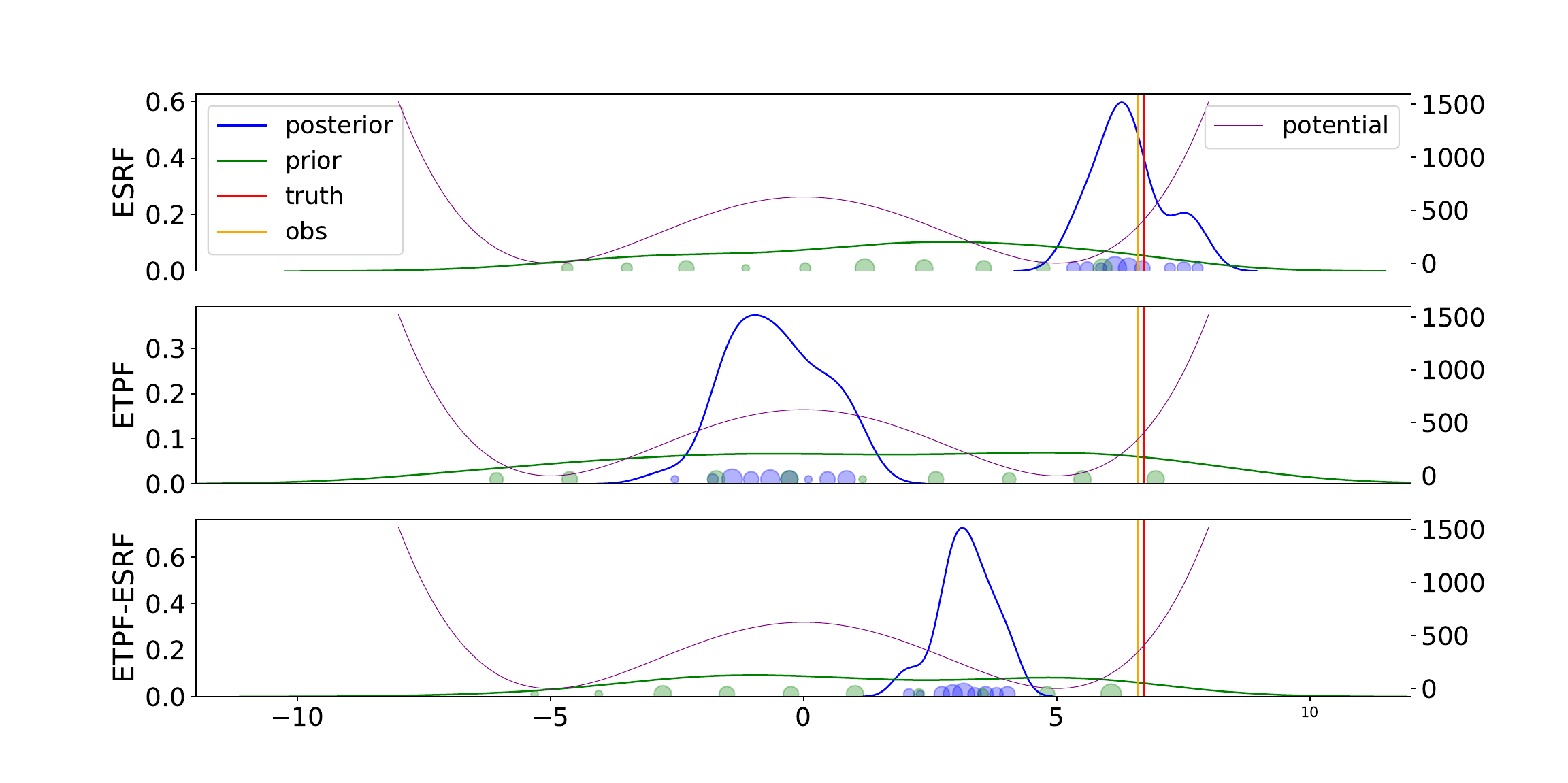}
\caption{Visualisation of the prior and posterior distributions as well as the truth and observation over two time steps (from left to right), highlighting the notable tail switching dynamics between time steps.}
\label{fig:Langevin}
\end{figure}

As the first prototypical scenario, we consider the evolution of the dampened Langevin equation:
\begin{align}
\begin{split}
    q_t &= M^{-1} p_t \, dt,\\
    p_t &= \nabla \phi(q_t) \, dt - \gamma p_t + \sigma M^{\frac{1}{2}}  dW_t.
\end{split}
\end{align}
The deterministic forces are given by the gradient of the double well potential \(\phi(z) = (z-5)^2(z+5)^2\), and the stochastic forces by the standard normal Wiener process \(W_t\).
We obtain the time discrete numerical approximation by the BAOAB splitting strategy as presented in \cite{Leimkuhler2012}. As stepsize we choose \(\Delta t = 0.01 \) and for the physical parameters we choose \(M=1\), \(\gamma=10\), \(\sigma = \sqrt{5000}\). Further, we observe all variables in observation intervals of $\Delta t_{obs}=0.8$ with observation error covariance $R=diag(0.5,0.5)$. For our purpose, a total of $50$ assimilation steps are performed. We include ensemble inflation and particle rejuvenation, see \Cref{sec:inflation,sec:rejunvenation}. We set the inflation parameter $\gamma=1.05$, the rejuvenation parameter $\tau=0.2$. In this setup, we compare ESRF, ETPF, Bootstrap, FPF and their combination controlled by the tempering criteria discussed in \Cref{sec:tempering}. \\
To understand the intricacies of the filter behaviour, we consider a concrete scenario: Consider a system of particles with mean exit time from a potential well comparable to the length of the observation interval. To investigate the performance of the filters and their tempered versions under a change to another potential well we initialise an ensemble of $N_{ens}=35$ particles.  
\Cref{fig:Langevin} provide a visual representation of the comparative performance between the ESRF, ETPF, and their tempered combination in the described scenario. Remarkably, the data indicates a notable advantage of the ESRF over the ETPF, particularly evident in situations where particles undergo tail switching dynamics, i.e. \Cref{fig:Langevin}. Further, in \Cref{fig:Langevin} we observe the tempered filter (ETPF-ESRF) to be competitive with the ESRF more likely than the single ETPF.
Here, we apply the three criteria for the tempering approach as described in \Cref{sec:tempering}. The threshold for the ESS criterion is set to $0.5$ and we choose a factor of $1.5$ for the IQR criterion. All numerical results are collected in \Cref{tab:Stat_Langevin_35}  \\

\begin{table}[htp]
    \centering
    \scalebox{0.9}{
    \begin{tabular}{|l|l|}
    \hline
        \text{Filter} &  \text{RMSE}   \\ \hline
         ESRF & \color{cherryred}{0.41761}  \\
         ETPF & \color{darkblue}{1.99017 } \\
         Bootstrap & \color{darkgreen}{1.25118}   \\ 
         FPF & \color{darkorange}{5.03897} \\ \hline
         
         ETPF-ESRF & 0.46150 \\
         Bootstrap-ESRF & 0.41130 \\ 
         FPF-ESRF &  0.41151   \\ \hline
         
         ESS-ETPF-ESRF & 0.46150  \\
         ESS-Bootstrap-ESRF & 0.41130  \\
         ESS-FPF-ESRF & 0.41151  \\ \hline

         IQR-ETPF-ESRF & \color{darkblue}{0.46150}\\
         IQR-Bootstrap-ESRF & \color{darkgreen}{0.41130}  \\
         IQR-FPF-ESRF & \color{darkorange}{0.41151}   \\ \hline
    \end{tabular}
    }
    \caption{Root Mean Squared Error (RMSE) for all filter combinations in the Langevin experiment.}
    \label{tab:Stat_Langevin_35}
\end{table}

\newpage
\subsection{Lorenz 63}\label{res:L63}
\begin{figure}[!ht]
\centering 
\includegraphics[width=1\textwidth]{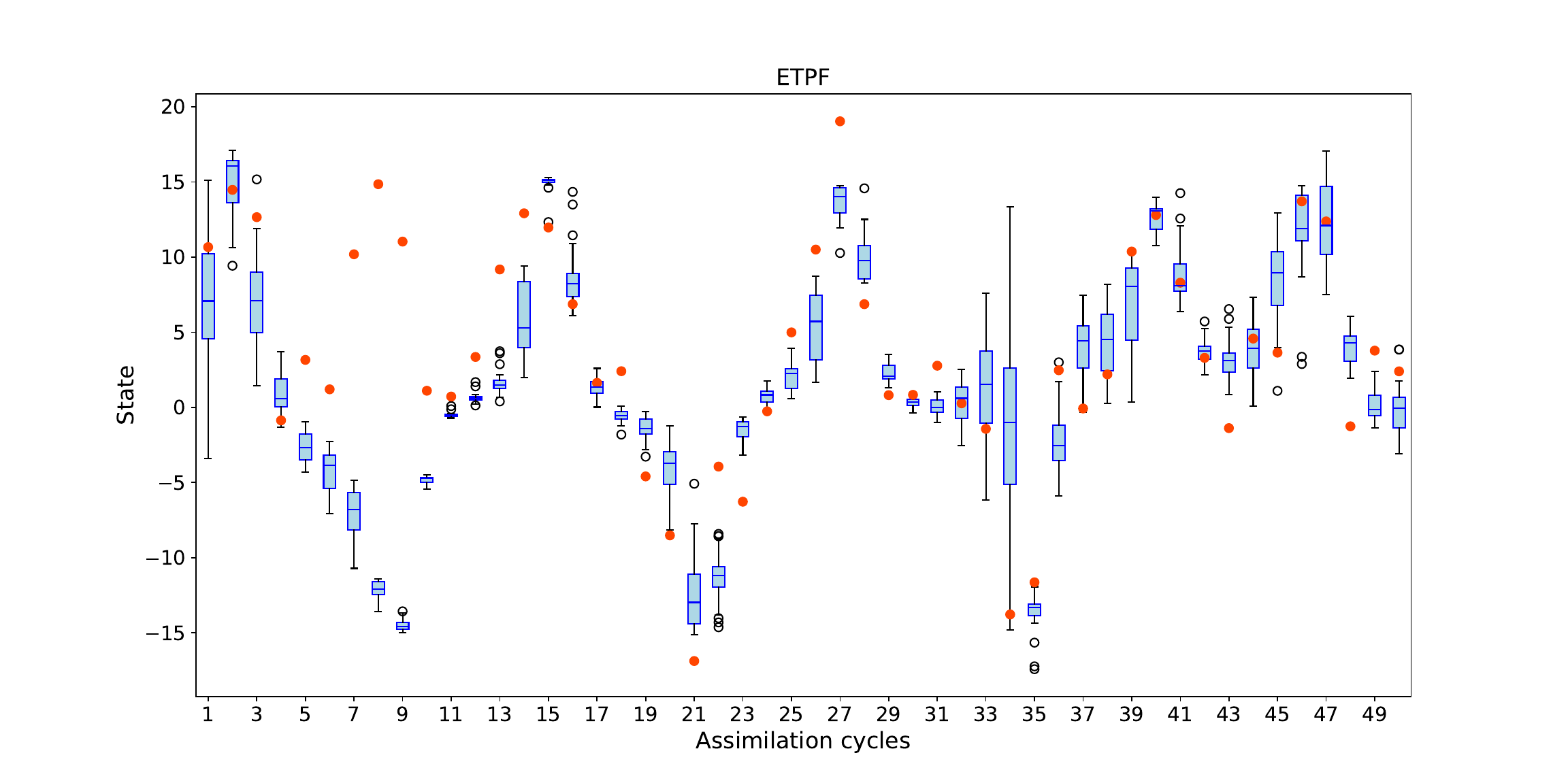}
\includegraphics[width=1\textwidth]{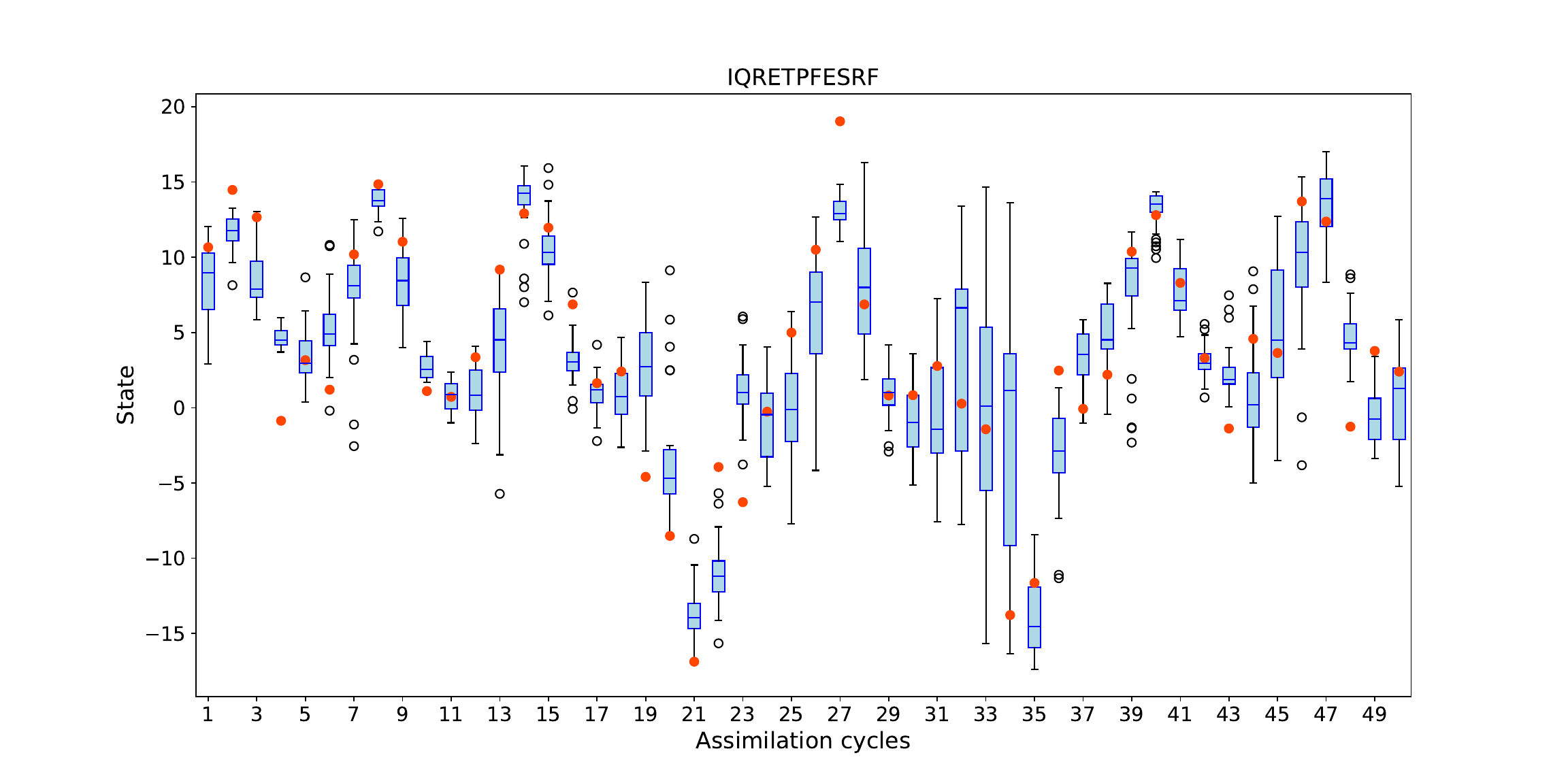}
\caption{The box plots include the observations and the ensemble state at the first state of the Lorenz 63 over a range of 50 assimilating cycles. 
The upper and lower diagram depicts the statistics computed by the plain ETPF and the tempered (IQR-criterion) version respectively.}
\label{fig:Boxplot_ETPFESRF_L63}
\end{figure}

We use the chaotic Lorenz 63 system with the standard parameter setting, whereby $\sigma=10, \rho=28,$ and $\beta=8/3$. We observe the first component of the three-dimensional system in observation intervals of $\Delta t_{obs}=0.12$ with observation error variance $R=8$. Our experiment covers $50000$ assimilation steps, whereby the first $500$ are neglected for the evaluation, to ensure the particles to be sufficiently close to the attractor. We include ensemble inflation and particle rejuvenation, see \Cref{sec:inflation,sec:rejunvenation}. Note that $\gamma=1.05$, $\tau=0.2$ as in the previous example. This data assimilation setting has already been used in \cite{Acevedo2017}. 
As in the previous example, we test several combinations of filters and 
present the results in \Cref{tab:Stats_L63_25}. Remarkably, the approach using the IQR performs best when comparing the RMSE of the single ETPF with the tempered version and beats even the single ESRF. In general, the IQR approach holds better results then the other two approaches.
\Cref{fig:Boxplot_ETPFESRF_L63} depicts first the box plots of the prior observations for the single ETPF. As visible in assimilation cycle $5-13$, the particles do not cover the observations. The second graph in \Cref{fig:Boxplot_ETPFESRF_L63} shows the tempered filter using the IQR criteria with quantile factor $1.5$. 

\begin{table}[!ht]
    \centering
    \scalebox{0.9}{
    \begin{tabular}{|l|l|l|l|l|l|}
    \hline
        \text{Filter} &  \text{RMSE}    \\ \hline
         ESRF & \color{cherryred}{2.10011}  \\
         ETPF & \color{darkblue}{3.55604}  \\
         Bootstrap & \color{darkgreen}{5.94889}  \\ 
         FPF & \color{darkorange}{4.94511}  \\ \hline

         ETPF-ESRF & 2.06520   \\
         Bootstrap-ESRF &2.08702    \\
         FPF-ESRF & 2.22518   \\ \hline
         
         ESS-ETPF-ESRF & 1.75024  \\
         ESS-Bootstrap-ESRF & 3.12627     \\
         ESS-FPF-ESRF & 2.22518   \\\hline

         IQR-ETPF-ESRF & \color{darkblue}{1.64179}     \\
         IQR-Bootstrap-ESRF &  \color{darkgreen}{2.01076}     \\
         IQR-FPF-ESRF & \color{darkorange}{2.22518}       \\\hline
    \end{tabular}
    }
    \caption{Root Mean Squared Error (RMSE) for all filter combinations in the Lorenz 63 model.}
    \label{tab:Stats_L63_25}
    
\end{table}

\newpage
\subsection{Lorenz 96}\label{res:L96}
\begin{figure}[!ht]
\centering 

\includegraphics[width=1\textwidth]{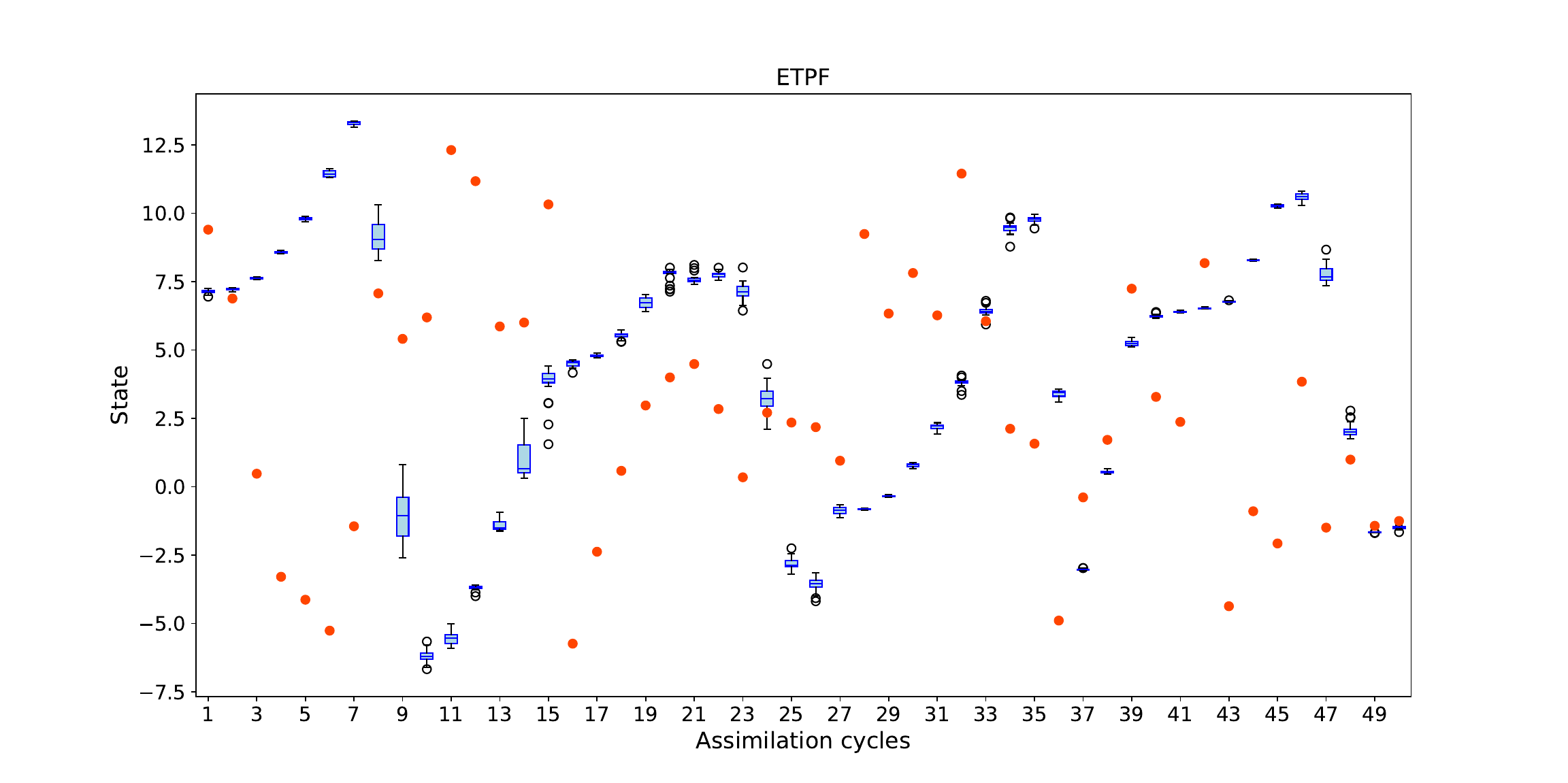}
\includegraphics[width=1\textwidth]{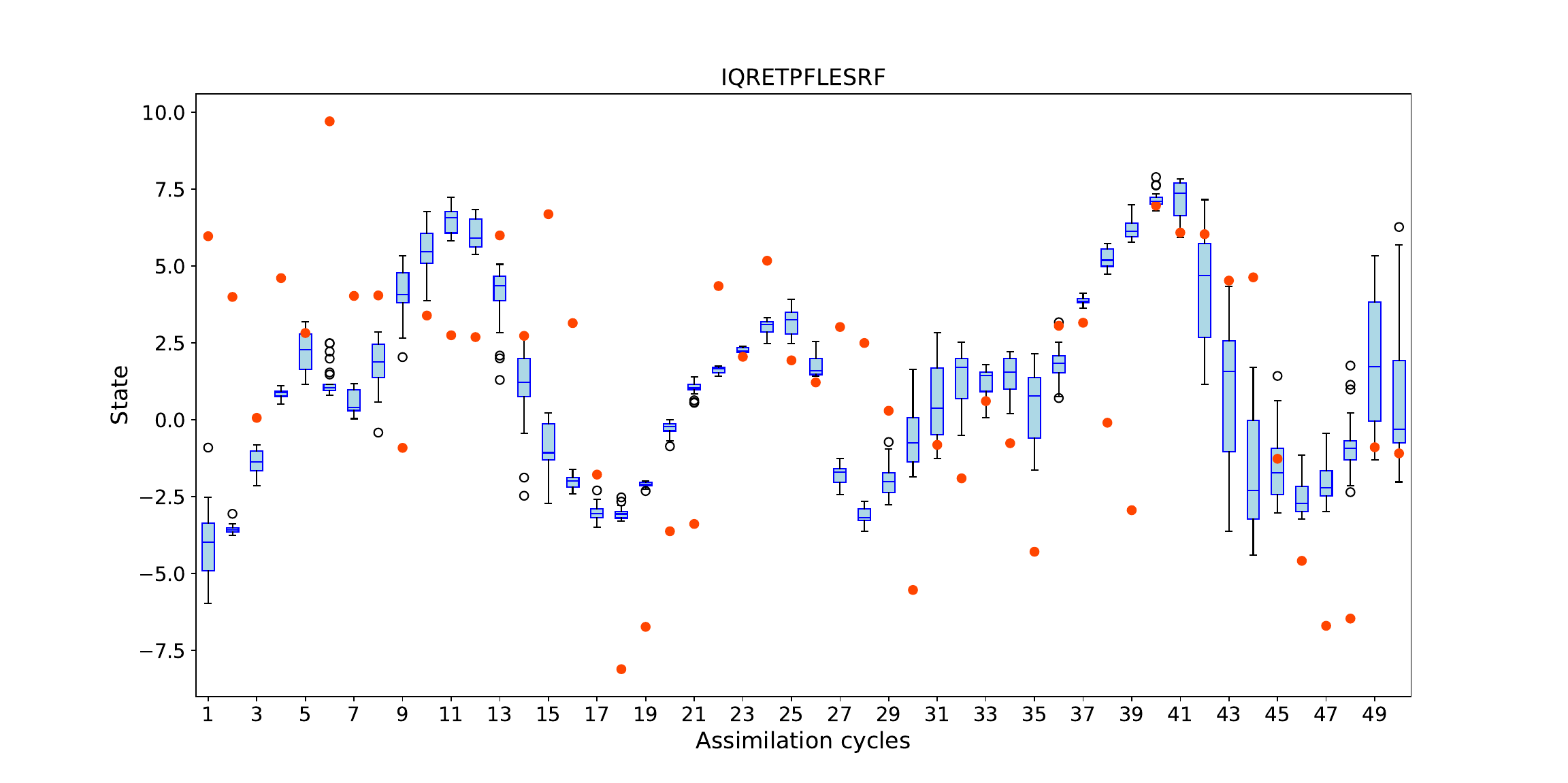}
\caption{The box plots include the observations and the ensemble state at the first grid point of the Lorenz 96 discretization over a range of 50 assimilating cycles. 
The upper and lower diagram depicts the statistics computed by the plain ETPF and the tempered (IQR-criterion) version respectively.}
\label{fig:Boxplot_ETPFESRF_L96}

\end{figure}

The discrete dynamical system introduced in the seminal work by \cite{Lorenz1996} is a widely adopted model to provide spatio-temporal benchmarks for data assimilation methods. The model originates from an unstable discretization of the dampened Burgers' equation with constant forcing. In contrast to a stable discretization, the evolution in time becomes chaotic.
We choose $120$ grid points, the forcing to be $8$ and observe every second grid point in observation intervals of $\Delta t_{obs}=0.11$. The observation error covariance is $R=diag(8,\dots,8)$ and the localisation radius is equal to two grid points. This experiment cover 5000 assimilation steps, whereby the first 50 are neglected. Further, $\gamma=1.05$, $\tau=0.2$, and $\alpha=0.2$.  
Inflation, rejuvenation, and localisation can also be applied in the tempered filters according to \Cref{sec:inflation,sec:rejunvenation,sec:localisation}. Comparing the localised filters, the RMSE varies more strongly. 
The minimum RMSE is obtained with the localised ETPF (LETPF) as visible in \Cref{tab:Stat_L96_35}. 
In general the tempered versions show a slightly better RMSE than the single filters alone.
Focusing on the ETPF in \Cref{tab:Stat_L96_35}, marked in blue, also shows that the different choices of tempering lead to similar values of the RMSE. 
This holds true for all other combinations as well and suggests the following heuristic. 
If the RMSE of the single filters differs too much or too less, then their tempered version is unable to substantially reduce the RMSE. 
Therefore, we refrain from combinations of the LETPF and LESRF, as they both perform equally well individually. 
We emphasise the substantial improvement for the FPF, marked in orange in \Cref{tab:Stat_L96_35}, especially in contrast to the other particle filters and their tempered versions. 
Finally, we observe again the increased spread in the ensemble in \Cref{fig:Boxplot_ETPFESRF_L96} for the plain ETPF and the tempered version using the LESRF and the IQR criteria with quantile factor zero. 
\\

\begin{table}[!ht]
    \centering
    \scalebox{0.9}{
    \begin{tabular}{|l|l|}
    \hline
        \text{Filter} &  \text{RMSE}    \\ \hline
         ESRF & 2.48515  \\
         LESRF & 1.08808  \\ 
         ETPF & \color{darkblue}{2.96073}   \\
         LETPF & 1.05996  \\ 
         Bootstrap & \color{darkgreen}{2.96868}   \\ 
         FPF &\color{darkorange}{2.15463}  \\ \hline

         ETPF-ESRF &2.86639   \\
         Bootstrap-ESRF &2.88486    \\
         FPF-ESRF & 1.98065   \\ \hline
         
         ESS-ETPF-ESRF & 2.86633  \\
         ESS-Bootstrap-ESRF & 2.91596   \\
         ESS-FPF-ESRF &  1.98306 \\ \hline

         IQR-ETPF-ESRF & 2.86639  \\
         IQR-Bootstrap-ESRF & 2.88486   \\
         IQR-FPF-ESRF &  1.98065   \\ \hline

         ETPF-LESRF & 2.51371\\
         Bootstrap-LESRF &2.60774  \\
         FPF-LESRF &  1.37967  \\ \hline
         
         ESS-ETPF-LESRF & 2.53425 \\
         ESS-Bootstrap-LESRF & 2.62288   \\
         ESS-FPF-LESRF &  1.37967  \\ \hline

         IQR-ETPF-LESRF & \color{darkblue}{2.51371}  \\
         IQR-Bootstrap-LESRF &\color{darkgreen}{2.60774}    \\
         IQR-FPF-LESRF &  \color{darkorange}{1.37967}  \\ \hline
    \end{tabular}
    }
    \caption{Root Mean Squared Error (RMSE) for all filter combinations in the Lorenz 96 model.}
    \label{tab:Stat_L96_35}
\end{table}

\newpage
\subsection{Shallow water}\label{res:SHW}
The shallow water system is a set of nonlinear hyperbolic partial differential equations that describe the fluid flow below a pressure surface. The equations are derived from depth-integrating the Navier-Stokes equations. Hence, conservation of mass implies that the vertical velocity scale of the fluid is small compared to the horizontal velocity scale. 
The characteristic speed of the shallow water waves depends on the water height above ground and which in turn leads to nonlinear effects e.g. wave breaking. 

Let $\Omega = [0,L] \subset \mathbb{R}$ then shallow water system with gravity $g$ and orography $z$ reads 
\begin{align}
\begin{split}
    h_t + \nabla \cdot hu&=0, \\
    (hu)_t + \nabla \cdot (hu \otimes u ) + gh\nabla h &= -gh \nabla z, 
\end{split}
\end{align}
\noindent 
where $h$ denotes the height over ground i.e. the total height is given by $H=h+z$ and $u$ denotes the two-dimensional velocity vector. On the boundaries, we assume a wall at \(x=L\) and homogenous Neumann boundary conditions at \(x=0\). A solution to the initial value problem, satisfying $u=\Vec{0}$ and $\nabla (z+h)= \Vec{0}$ for all times is called the lake at rest steady state.

Despite being reasonably popular for the modelling of nonlinear water waves and internal gravity or Rossby waves in the atmosphere \cite{Majda2003}, the shallow water equations, a-priori, are a questionable choice as a model for data assimilation due to the lack of any dissipative mechanism as e.g. in the Navier-Stokes case (c.f. \cite{Jones1992}). 
Nevertheless, the model is sufficient to serve as a high-dimensional example for strong nonlinear effects. In contrast to the Lorenz 96 model, the spatial solutions are stable and allow for a meaningful physical interpretation.

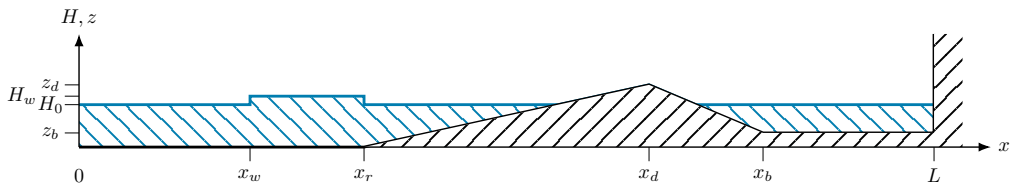
\begin{figure}[ht]
    \centering
    \scalebox{0.75}{
    \begin{tikzpicture}
        \draw[darkblue, line width=0.5mm] (0,0.75)-- (3,0.75) -- (3,0.9) -- (5,0.9) -- (5,0.75) -- (8.41,0.75) -- (10,1.1) -- (10.82,0.75) -- (15,0.75);
        \fill[pattern={Lines[
                  distance=2mm,
                  angle=-45,
                  line width=0.25mm
                 ]}, pattern color=darkblue] (0,0)--(0,0.75)-- (3,0.75) -- (3,0.9) -- (5,0.9) -- (5,0.75) -- (8.41,0.75) -- (10,1.1) -- (10.82,0.75) -- (15,0.75)--(15,0);
 
        \draw[black, line width=0.5mm] (0,0) -- (5,0) -- (10,1.1) -- (12,0.25) -- (15,0.25) -- (15,2);
        \fill[white] (0,0) -- (5,0) -- (10,1.1) -- (12,0.25) -- (15,0.25)--(15,2)--(15.5,2)--(15.5,0);
        \fill[pattern={Lines[
                  distance=2mm,
                  angle=45,
                  line width=0.25mm
                 ]}, pattern color=black] (0,0) -- (5,0) -- (10,1.1) -- (12,0.25) -- (15,0.25)--(15,2)--(15.5,2)--(15.5,0);
        \draw[-Latex,line width=0.25mm] (0,0) -- (16,0) node[right]{$x$};
        \draw[-Latex,line width=0.25mm] (0,0) -- (0,2) node[above]{$H,z$};
        \draw (0.0,-0.5) node{$0$};
        \draw (3,0)  --  (3,-0.25);
        \draw (3,-0.50)  node{$x_w$};
        \draw (5,0)  --  (5,-0.25);
        \draw (5,-0.50)  node{$x_r$};
        \draw (10,0) -- (10,-0.25);
        \draw (10,-0.50) node{$x_d$};
        \draw (12,0) -- (12,-0.25);
        \draw (12,-0.50) node{$x_b$};
        \draw (15,0) -- (15,-0.25);
        \draw (15,-0.50) node{$L$};
        \draw (0,1.1) -- (-0.25,1.1);
        \draw (-0.5,1.1) node{$z_d$};
        \draw (0,0.75) -- (-0.25,0.75);
        \draw (-0.5,0.75) node{$H_0$};
        \draw (0,0.9) -- (-0.25,0.9);
        \draw (-1,0.9) node{$H_w$};
        \draw (0,0.25) -- (-0.25,0.25);
        \draw (-0.5,0.25) node{$z_b$};
    \end{tikzpicture}
    }
    \caption{Initial condition of the used shallow water system. }
    \label{fig:shw_ic}
\end{figure}
\begin{figure}[ht]
\centering 
\includegraphics[width=0.9\textwidth]{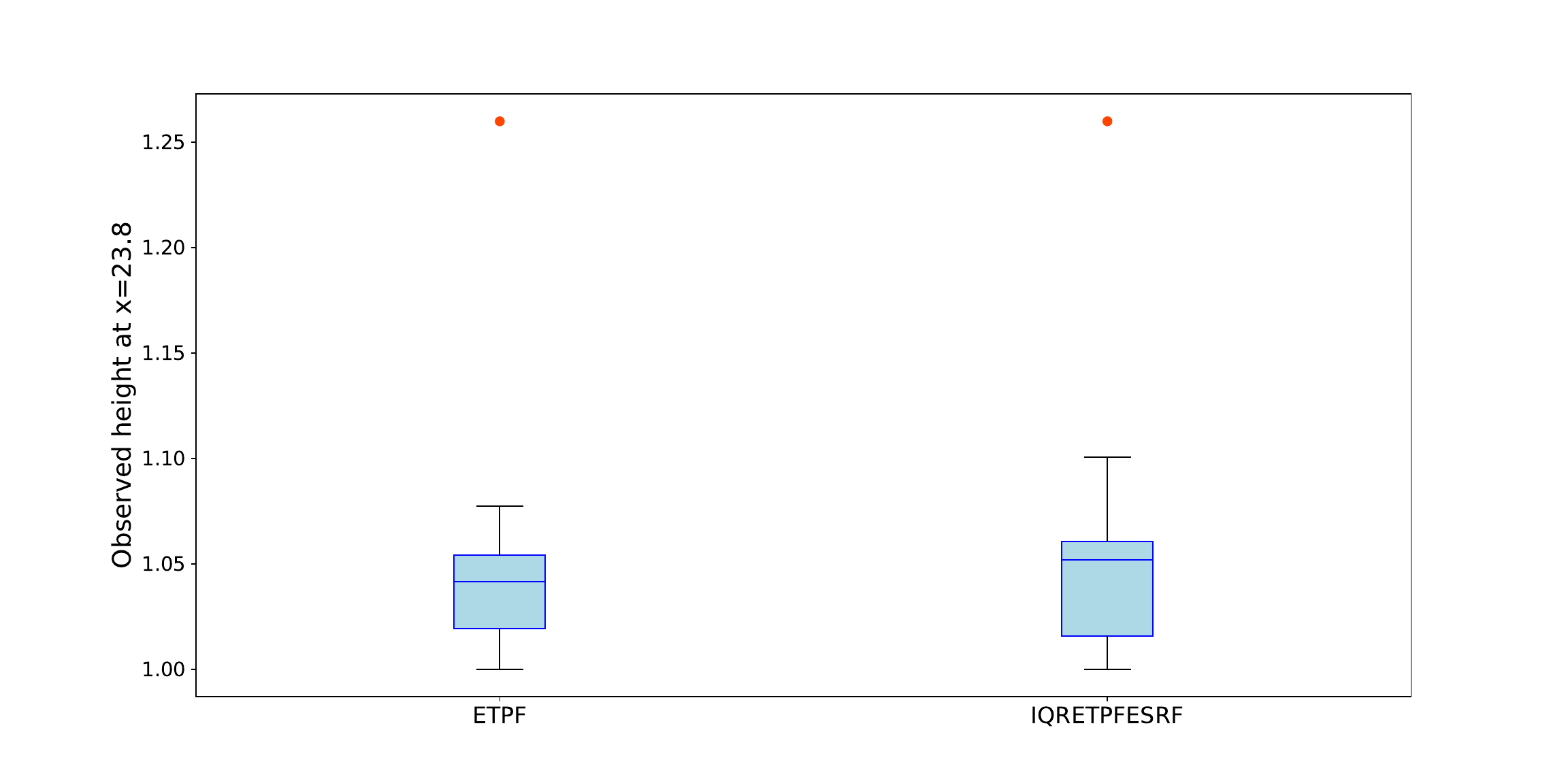}
\caption{The box plots include the observations and the ensemble state at the height $x=23.8$ of the shallow water system comparing the ETPF and the IQR version respectively. }
\label{fig:Boxplot_SHW}
\end{figure}

We compute a numerical approximation of solutions to the initial valued problem of the shallow water equations by a second-order centred finite volume scheme presented in \cite{Bouchut2007} using hydrostatic reconstruction developed in  \cite{Audusse2004,Botta2004}.
This enables us to sample from the prior distribution efficiently and robustly, even in the given scenario with wet-dry interfaces, i.e. when the non-negative \(h\) vanishes in some parts of the domain. Furthermore, the finite volume method preserves the mass of the fluid exactly.
It is important to denote that although the actual filtering distribution is non-negative, some filter modifications as rejuvenation might set samples of the posterior to negative values. In this case, we cut off the values at zero. 
We choose the initial condition of the reference solution according to \Cref{fig:shw_ic}. The discretization is spatially uniform on \([0,L-x_b)\) with $49$ grid cells and contains a single grid cell covering \( [L-x_b,L]\). The specific values are \(x_r=10\), \(x_d=20\), \(x_b=24\), \(L=28\), \(H_0=1.5\), \(H_w=1.8\), \(z_b=0.5\), \(z_d = 2.2\). For the reference solution \(x_w=3\) and the ensemble height is distributed according to \(x_w \approx \mathcal{N}(6,1)\), i.e. the initial condition is misspecified on the \(99.87\%\) percentile.   
The velocity of the initial ensemble is generated by 
\(v = \sqrt{g(H_w-H_0)}\), where \(g=9.81\).
The observational window is given by \(\Delta t_{obs}=2.5\) and we observe the fluid height at the centre of every 8\textsuperscript{th} grid cell. We run 3 assimilation cylces, but skip the first for collecting the statistics. Note that $\gamma=1.05$, $\tau=0.2$, and $\alpha=0.2$ as usual.
We evaluate the RMSE in the large cell covering the basin behind the dam and average over the two assimilation cycles. The collected results are listed in \Cref{tab:Stat_SHW}.
Furthermore, we depict the ensemble statistics for the rightmost point of observation in \Cref{fig:Boxplot_SHW}. Although neither the RMSE nor the ensemble statistics change drastically, we observe a similar trend as in previous examples. Interestingly enough, we can spatially locate the differences between the plain and the tempered version of the filter, where the particle filters, by construction of the orography, have reduced spread in the distribution. 

\begin{table}[!ht]
    \centering
    \scalebox{0.9}{
    \begin{tabular}{|l|l|}
    \hline
        \text{Filter} &  \text{RMSE}   \\ \hline
         ESRF &  \color{cherryred}{0.038228}    \\
         ETPF & \color{darkblue}{0.049417}    \\
         Bootstrap & \color{darkgreen}{0.056183}     \\ 
         FPF & \color{darkorange}{0.071165}     \\ 
         LESRF & 0.048350   \\
         LETPF & 0.056134    \\ \hline

         ETPF-ESRF & \color{darkblue}{0.040833} \\
         Bootstrap-ESRF & 0.030875 \\
         FPF-ESRF & 0.053365   \\ \hline
         
         ESS-ETPF-ESRF & \color{darkblue}{0.051588}  \\
         ESS-Bootstrap-ESRF & 0.056429  \\
         ESS-FPF-ESRF & 0.053365   \\ \hline

         IQR-ETPF-ESRF & \color{darkblue}{0.040833}  \\
         IQR-Bootstrap-ESRF & \color{darkgreen}{0.030875}  \\
         IQR-FPF-ESRF & \color{darkorange}{0.053365}  \\ \hline
    \end{tabular}
    }
    \caption{Root Mean Squared Error (RMSE) for all filter combinations in the shallow water experiment.}
    \label{tab:Stat_SHW}
\end{table}

\newpage
\section{Discussion}
The focus of this study is to better the performance of particle filters, especially those described in \Cref{sec:Bootstrap,sec:ETPF,sec:FPF}, by the means of tempering \Cref{sec:tempering}. The novel approach presented in \Cref{sec:IQR}, receives the greatest attention and differs from the approaches developed in \cite{Frei2013,Chustagulprom2016,Nerger2021}.   
Specifically, \cite{Nerger2021} posed two choices of the bridging parameter for the tempering approach, First, depending on the effective sample size as recommended in \cite{Frei2013},
second, on skewness and kurtosis of the empirical distribution in observational space. This should quantify the non-Gaussianity of the ensemble distribution and as pointed out in, this shows promising results, but requires further fine-tuning. 
This work is not primarily concerned with the bridging parameter, but rather with the question whether to temper or not. 
Targeted to the challenging situations for particle type filters, we proposed a novel tempering criterion, based on descriptive statistics, which is straightforward to compute and has a clear interpretation, i.e. if the observation is within the interquartile range of the observed particles.
This criterion is complementary to a threshold on the effective sample size, as the latter does not take the observation into account.
As particle filters in high-dimensional systems usually suffer from the lack of overlapping support of the likelihood and the empirical prior distribution, we propose to temper the likelihood by a Gaussian-approximative filter to push the particles possibly towards the observation. 
This prototypical situation is illustrated in our experiment by applying several particle filters to a filtering scenario based on the double well potential c.f. \Cref{fig:Langevin} and \Cref{tab:Stat_Langevin_35}. As \Cref{fig:Langevin} visualises, the ETPF struggles clearly to switch the modes of the double well. Furthermore, it is visible that the tempered approach shifts the ensembles closer to the observation.
\Cref{tab:Stat_Langevin_35} points out significant improvements for all three particle type filters, highlighted in blue for the ETPF, green for the Bootstrap and orange for the FPF. Note that there is no indication which of the tempering approaches serves as best one. 
Further, we test the tempering approaches on different particle filters, tuning them with rejuvenation and localisation as needed, and provided tests on the defacto standard benchmarks based on both Lorenz models.  
In the case of the Lorenz 63 model, we observe improvements for the time-averaged root mean square error, for all tempered versions, but find the IQR criterion to perform best as listed in \Cref{tab:Stats_L63_25}. IQR-ETPF-ESRF demonstrated the lowest RMSE ($1.64179$), surpassing even the performance of the single ESRF ($2.10011$). The IQR criteria also outperform the other two tempering approaches.  This finding counts also for the other two particle type filters as highlighted in \Cref{tab:Stats_L63_25}. 
In the Lorenz 96 experiments, the improvements are clearly visible, but not as pronounced as in the previously mentioned case.  Considering \Cref{tab:Stat_L96_35}, the LETPF exhibited the lowest RMSE ($1.05996$), demonstrating that localisation approaches can increase filter accuracy. 
As previously mentioned, we focus on the particle type filters (without localisation) and their tempered variations. \Cref{tab:Stat_L96_35} shows very small improvements for the ETPF and Bootstrap. 
Comparing the RMSE of the IQR-FPF-LESRF ($1.37967$), which outperforms the FPF ($2.15463$), exhibited a substantial improvement in performance, suggesting that the incorporation of IQR alongside FPF and LESRF can yield significant enhancements. 
Further, \Cref{fig:Boxplot_ETPFESRF_L96} illustrates the advantages of using IQR, highlighting its effectiveness in enhancing filtering accuracy.
Finally, we provide a test case based on the shallow water equations with orography. Here, we investigate essentially two cycles of a constructed filtering scenario and would like to highlight, that this scenario provides a situation, where the criterion based on the ESS provides worse results than the constant tempering or the IQR based criterion, i.e. \Cref{tab:Stat_SHW}. 
This points out the robustness of the IQR criterion across all experiments in contrast to the ESS, which does not target the specific difficulties under investigation.
Furthermore, \Cref{fig:Boxplot_SHW} illustrates again the shift towards the observation, localised at the specific subspace (gridpoint).

\section{Conclusion}
Tempering significantly enhances the performance of filtering methods, particularly for approximate consistent filters. Additionally, it can help stabilise numerically unstable filters through occasional robust updates, using a combination of filters rather than a single fixed choice. However, determining the optimal timing for employing these tools during filtering procedures remains unclear. We propose a novel approach to adaptively decide on the use of an approximate tempering step, using the interquartile range (IQR). To numerically validate this tempering criterion across various filters, we consider three representative filters from key families of consistent filters: standard sequential Monte Carlo resampling, Transport-map-induced, and McKean-Vlasov filters. Each of these filters is paired with an established representative from the class of Gaussian approximate filters. Overall, the new adaptive tempering approach of the considered consistent filters paired with adaptive tempering steps resulted in lower or equal RMSE compared to the performance of the individual filter, a fixed tempering cycle, and an additional effective sample size criterion. This effect is evident across three different types of toy models, each with a specific focus. The promising results obtained in this study suggest further investigation into more realistic filtering scenarios based on atmospheric fluid dynamics. Another important aspect will be a careful examination of the numerical performance of the Feedback Particle Filter in the non-asymptotic regime with small ensemble sizes.

\section*{Acknowledgments}
This research was funded in whole or in part by the Austrian Science Fund (FWF) [10.55776/DOC78]. For open access purposes, the author has applied a CC BY public copyright license to any author-accepted manuscript version arising from this submission.   \\
The research of Gottfried Hastermann, Iris Rammelmüller and Jana de Wiljes has been partially funded by the Deutsche Forschungsgemeinschaft (DFG)- Project-ID 318763901 - SFB1294. \\
The program code is available as ancillary file from the arXiv page of this paper.

\bibliographystyle{plainurl}
\bibliography{library}

\section{Appendix}
\noindent
\begin{algorithm}[]
\small
    \caption{Ensemble Square Root Filter (ESRF), see also the presentation in \cite{Reich2015} and \Cref{sec:ESRF}}  \label{alg:ESRF}
   \begin{algorithmic}[1]
        \vspace{0.25\baselineskip}
        \Require{ensemble forecast $z^f_i\in \mathbb{R}^{N_x}, \quad  \forall i\in \{1, \dots,N_{ens}\}$ and observation $y^{obs} \in \mathbb{R}^{N_y}$, $\gamma=1.05$}
        \vspace{0.5\baselineskip} 
        \State $\overline{z}^f =  \frac{1}{N_{ens}} \sum_{i=1}^{N_{ens}} z^f_{i}$
        \vspace{0.5\baselineskip}
        \State $A^f = (z^f_{i} -\overline{z}^f)$
        \vspace{0.5\baselineskip}
         \State inflation: $z_i^f = \overline{z}^f + \gamma A^f$
        \vspace{0.5\baselineskip}
        \State compute step 1 and 2 again
        \vspace{0.5\baselineskip}
        \State $S = \left( I_{N_{ens}} + \frac{1}{N_{ens}-1} (HA^f)^T R^{-1}  (HA^f)  \right)^{-1/2} $
        \vspace{0.5\baselineskip}
        \State $ \hat{w}_i =\frac{1}{N_{ens}}- \frac{1}{N_{ens}-1}  e_i^T S^2 (HA^f)^T R^{-1} (H\overline{z}^f-y^{obs}) $
        \vspace{0.5\baselineskip}
        \State $ d_{ij} = \hat{w}_i - \frac{1}{N_{ens}} + s_{ij} $
        \vspace{0.5\baselineskip}
        \State $z_{j}^{a} =\sum_{i=1}^{N_{ens}}  z^f_{i} d_{ij} $\\
        \vspace{0.75\baselineskip}
        \Return $z^a_{j}$ for every $j\in \{1, \dots,N_{ens}\}$
   \end{algorithmic} 
\end{algorithm}

\begin{algorithm}[ ]
\small
    \caption{Bootstrap Particle Filter, see \Cref{sec:Bootstrap}} \label{alg:Bootstrap}
   \begin{algorithmic}[1]
        \vspace{0.25\baselineskip}
        \Require{initial ensemble $z^f_{i}\in \mathbb{R}^{N_x}, \quad  \forall i\in \{1, \dots,N_{ens}\}$ and observation $y^{obs} \in \mathbb{R}^{N_y}, \tau=0.2$}
        \vspace{0.5\baselineskip} 
        \State $\overline{z}^f =  \frac{1}{N_{ens}} \sum_{i=1}^{N_{ens}} z^f_{i}$
        \vspace{0.5\baselineskip}
        \State $A^f = (z^f_{i} -\overline{z}^f)$
        \vspace{0.5\baselineskip}
        \State $ \hat{w}_{i}\propto \exp \left( - \frac{1}{2}(Hz^f_{i}- y^{obs})^\top  R^{-1} (Hz^f_{i}- y^{obs}) \right) \cdot w_0, \ \ \ w_0=\frac{1}{N_{ens}} $
        \vspace{0.5\baselineskip}
        \State $w = \frac{\hat{w}_{i}}{\sum_{i=1}^{N_{ens}} \hat{w}_{i}}$
        \vspace{0.5\baselineskip}
        \State $rejuvenation= \frac{\tau }{N_{ens}-1} A^f \xi, \ \xi \sim \mathcal{N}(0,1)$
        \vspace{0.5\baselineskip}
        \State $z^a_{i}= rng.choice(z^f_{i}, N_{ens}, p=w, axis=1) + rejuvenation $\\
        \vspace{0.75\baselineskip}
        \Return $z^a_{i}$ for every $i\in \{1, \dots,N_{ens}\}$
   \end{algorithmic} 
\end{algorithm}

\begin{algorithm}[]
\small
\caption{Ensemble Transform Particle Filter (ETPF), see also the presentation in \cite{Reich2015} and \Cref{sec:ETPF}} \label{alg:ETPF}
\begin{algorithmic}[1]
\vspace{0.25\baselineskip}
\Require{ensemble forecast $z^f_i\in \mathbb{R}^{N_x}, \quad  \forall i\in \{1, \dots,N_{ens}\}$ and observation $y^{obs} \in \mathbb{R}^{N_y}, \tau=0.2$}
\vspace{0.5\baselineskip}
\State $\overline{z}^f =  \frac{1}{N_{ens}} \sum_{i=1}^{N_{ens}} z^f_{i}
$
\vspace{0.5\baselineskip}
\State $A^f = (z^f_{i} -\overline{z}^f)$
\vspace{0.5\baselineskip}
\State  $ w_{i} \propto \exp \left( - \frac{1}{2}(Hz^f_{i}- y^{obs})^\top  R^{-1} (Hz^f_{i}- y^{obs}) \right) $ 
\vspace{0.5\baselineskip}
\State $\hat{w}_{i} = \frac{w_{i}}{\sum_{i=1}^{N_{ens}} w_{i}}$ 
\vspace{0.5\baselineskip}
\State $T^\ast \gets$ solution of \eqref{eq:optimal_transport} using the Euclidean metric i.e. \(d=d_2\).  
\vspace{0.5\baselineskip}
\State $rejuvenation= \frac{\tau }{N_{ens}-1} A^f \xi, \ \xi \sim \mathcal{N}(0,1)$
\vspace{0.5\baselineskip}
\State $z^a_{i} = \sum_{i=1}^{N_{ens}} z^f_{i} N_{ens}t_{ij}^\ast + rejuvenation$\\
\vspace{0.75\baselineskip}
\Return $z^a_{i}$ for every $i\in \{1, \dots,N_{ens}\}$
\end{algorithmic}
\end{algorithm}

\begin{algorithm}[]
\small
    \caption{Feedback Particle Filter (FPF), see also the presentation in \cite{Taghvaei2020} and \Cref{sec:FPF}} \label{alg:Feedback}
   \begin{algorithmic}[1]
        \vspace{0.25\baselineskip}
        \Require{ensemble forecast $z^f_i\in \mathbb{R}^{N_x}, \quad  \forall i\in \{1, \dots,N_{ens}\}$ and observation $y^{obs} \in \mathbb{R}^{N_y}$, $\epsilon=0.01$, and observation operator $H(\cdot)$, $\phi_0$, $\mathcal{T}$ iteration, $\tau=0.2$}
        \vspace{0.5\baselineskip} 
        \State $g_{ij} = \exp{\frac{- \vert z^f_{i} - z^f_{j} \vert^2}{4 \epsilon}}$
        \vspace{0.5\baselineskip}
        \State $l_i = \exp \left( - \frac{1}{2}(Hz^f_{i}- y^{obs})^\top  R^{-1} (Hz^f_{i}- y^{obs}) \right)$
        \vspace{0.5\baselineskip}
        \State $ k_{ij} = \frac{g_{ij}}{\sqrt{\sum_lg_{il}}\sqrt{\sum_lg_{lj}}} \ l_j $
        \vspace{0.5\baselineskip}
        \State $d_i = \sum_j k_{ij}$
        \vspace{0.5\baselineskip}
        \State $T_{ij} = \frac{k_{ij}}{d_i}$
        \vspace{0.5\baselineskip}
        \State $ \pi_i = \frac{d_i}{\sum_j d_j}$
        \vspace{0.5\baselineskip}
        \State $ \hat{H} = \sum_{i=1}^{N_{ens}} \pi_j H(z^f_{i})$
        \vspace{0.5\baselineskip}
        \State$\phi_t = \textnormal{solution of }  (\mathop{id} - T + \mathbf{1}) \phi_t = \varepsilon R^{-1} (H(z^f)-\hat{H})$
        \vspace{0.5\baselineskip}
        \State $r_i = \phi_\mathcal{T} + \epsilon H(z^f) $
        \vspace{0.5\baselineskip}
        \State $s_{ij} = \frac{1}{2 \epsilon} T_{ij} (r_j - \sum_{k=1}^{N_{ens}} T_{ik} r_k )$
        \vspace{0.5\baselineskip}
        \State $K_i = \sum_j s_{ij} z^f_{j} $
        \vspace{0.5\baselineskip}
        \State $A^f = (z^f_{i} -\overline{z}^f)$
        \vspace{0.5\baselineskip}
        \State $rejuvenation= \frac{\tau }{N_{ens}-1} A^f \xi, \ \xi \sim N(0,1)$
        \vspace{0.5\baselineskip}
        \State $z^a = z^f - K_i (\frac{1}{2}(H(z^f)-H(\overline{z_k}^f))-y^{obs}) + rejuvenation $ \\
        \vspace{0.75\baselineskip}
        \Return $z^a_{i}$ for every $i\in \{1, \dots,N_{ens}\}$
   \end{algorithmic} 
\end{algorithm}

\end{document}